\documentclass[10pt]{article}
\usepackage{amsmath,amsfonts,amsthm,amssymb,amscd}

\newcommand{\be}{\begin{equation}}
\newcommand{\ee}{\end{equation}}
\newcommand{\bs}{\begin{split}}
\newcommand{\es}{\end{split}}
\newcommand{\ba}{\begin{align}}
\newcommand{\ea}{\end{align}}
\newcommand{\basl}[1]{\begin{align}\begin{split}\label{#1}}
\newcommand{\bas}{\begin{align}\begin{split}}
\newcommand{\bl}[1]{\begin{equation}\label{#1}}
\newcommand{\el}{\end{equation}}

\newtheorem{theo}{Theorem}[section]
\newtheorem{prop}[theo]{Proposition}

\newtheorem{defi}[theo]{Definition}

\newcommand\N{\mathbb{N}}

\newcommand\R{\mathbb{R}}
\newcommand\C{\mathbb{C}}

\author{L. Amour\textsuperscript{1}, R. Lascar\textsuperscript{2}
and
J. Nourrigat\textsuperscript{1}}

\date{ \ \textsuperscript{1}Universit\'e of Reims\hskip 1cm \
\textsuperscript{2}Université Paris 7, Denis Diderot  }

\date{\ \textsuperscript{1}Universit\'e of Reims\hskip 1cm \
\textsuperscript{2}Universit\'e Paris 7, Denis Diderot  }

\title{Weyl calculus in Wiener spaces and in QED}

\begin{document}

\maketitle

\begin{abstract}
The concern of this article is a semiclassical Weyl calculus on an infinite dimensional Hilbert space $H$.
If $(i, H, B)$ is a Wiener triplet associated to $H$, the quantum state space  will be the space  of $L^2$
functions on $B$ with respect to  a Gaussian measure with $h/2$ variance, where $h$ is the
semiclassical parameter. We prove the boundedness  of our pseudodifferential operators (PDO) in the spirit
of Calder\'on-Vaillancourt with an explicit bound,  a Beals type characterization, and metaplectic covariance.
An application to a model of quantum electrodynamics (QED) is added in the last section (section \ref{s7}), for
fixed spin $1/2$ particles interacting with the quantized electromagnetic field (photons). We prove that
some observable time evolutions, the spin evolutions, the magnetic and  electric evolutions when subtracting their free evolutions, are PDO in our class.

\end{abstract}

\parindent=0pt

\



\tableofcontents

\parindent = 0 cm

\parskip 10pt
\baselineskip 14pt

\section{Introduction.}\label{s1}

Pseudodifferential operators theory, starting from 1965, has been adapted to Weyl quantization
by H\"ormander \cite{Hor} and Unterberger \cite{Unt-OH}. It is one of the main tools in  analysis of PDE
and spectral theory. See for instance H\"ormander \cite{Hor}, Lerner \cite{Ler},
Taylor \cite{Tay} and Zworski \cite{Zwo}.

In the works \cite{AJN-JFA}, \cite{ALN-BEALS} and \cite{ALN-QED-A}, an adaptation to  infinite
dimensional Hilbert spaces is proposed, with applications in  quantum electrodynamics (QED)
in \cite{ALN-QED-A}.

The goal of this article is to propose an adaptation of this infinite dimensional Weyl
calculus to a new class of symbols. The phase space is $H^2$, where $H$ is a separable Hilbert
space.  The   new class of symbols used here, denoted by $S(H^2, Q)$,  is defined in \cite{Jag}
using a given positive quadratic form $Q$ on $H^2$,  in spirit of H\"ormander \cite{Hor} or
Unterberger \cite{Unt-OH} works.  In this way
we obtain improvements with respect to \cite{AJN-JFA} and \cite{ALN-BEALS},  about calculus covariance and about applications to QED.

First we prove estimates of the $L^2$ norm of bounded operators
associated with functions $F$ in $S(H^2,Q)$ in terms of a Fredhlom
determinant. When the dimension is finite it makes more precise a classical
result of Calder\'on Vaillancourt \cite{CV} in the framework of H\"ormander
\cite{Hor} Weyl calculus. Our bounds are dimension free and depend only
in an explicit way of the quadratic form $Q$. Moreover a Beals type
characterization of PDO (see \cite{Bea}, \cite{Bon1}, \cite{Bon2}, \cite{BC}) is proved with as a consequence result about operators composition.

Another improvement with respect to \cite{AJN-JFA} is that our
class of symbols is coordinates free and in Section \ref{s4} we discuss the
covariance by the reduced symplectic group (see Shale \cite{Sha} and
B. Lascar \cite{Las}) at symbols level and by the metaplectic transforms
(named also Bogoliubov transforms) at operators one. See B. Lascar
\cite{Las} for $L^2$ framework or F. Hiroshima, I. Sasaki, H. Spohn, A.
Suzuki \cite{HSSS} and Bruneau Derezinski \cite{BD} for the abstract Fock
one.

The operators discussed here may be applied in QED theory if one
looks at interaction between the quantized filed and $N$ fixed particles
with $1/2$ spin in presence of a constant magnetic field. As seen in
Section \ref{s7}, particles spin, components of electric and magnetic
fields at each point $x\in\R^3$ and photon number may be viewed
through Weyl calculus. Such results were obtained in \cite{ALN-QED-A} with
an additional assumption on the infrared cut-off, no longer required
here, but symbols used  in \cite{ALN-QED-A} may  have
their supports confined in balls,  which is
not allowed here in view of the analytic properties of our symbols.

We define first the  class of symbols.

\begin{defi}\label{d1.1}(\cite{Jag}).   Let
$H$ be a real separable Hilbert, $Q$ a positive quadratic form on
$H^2$ and $S (H^2, Q) $ the class of functions $f \in C^{\infty }
(H^2)$  such that  there exists $C(f) >0$ satisfying for any integer $m\geq
0$, \be\label{1.1}|(d^m f ) (x) (U_1 , \dots , U_m ) | \leq C (f)
 Q( U_1) ^{1/2} \dots  Q( U_m) ^{1/2}.
 \ee
If $E$ is a Banach space we denote by $S (H^2, Q, E)$ the analogous
space of symbols taking values in $E$. The smallest constant satisfying
$(\ref{1.1})$ is denoted by $\Vert f \Vert _Q$. If ${\rm dim}\,H$ is
finite, the class $S_p (H^2, Q) $ of functions $f \in C^{p} (H^2)$
satisfying $(\ref{1.1})$ for all $m\leq p$ is also used and the best
constant $C(f)$ is denoted by $\Vert f \Vert _{p,Q}$.
\end{defi}

In what follows, $Q$ will be, \be\label{1.2}Q_A(x) = < Ax , x >, \ee
with $A\in {\cal L}(H^2)$, $A=A^*$, $A$  nonnegative trace class in
 $H^2$. We note that nonnegative quadratic form $Q$ can also be considered in the above definition.

Let us observe that the elements of $S(H^2, Q)$ are analytic
functions.
Let us recall that a class of operators with analytic symbols is introduced in  \cite{BMK}.
The Weyl operators associated with symbols of $S(H^2, Q)$
will be bounded operators in an $L^2$ space obtained by a below classic
result of Gross and Kuo.

\begin{theo}\label{t1.2}  (Gross \cite{Gro}, Kuo \cite{Kuo}).
For any real separable Hilbert space there exists a Banach space $B$
containing $H$, non unique, with a continuous injection
$i: H\hookrightarrow B$ with a dense range and  satisfying  $B'\hookrightarrow
H'=H \hookrightarrow B$, and for any $h>0$, there exists a probability
measure $\mu_{B , h}$ on the Borel $\sigma-$algebra of $B$
satisfying,
 \be\label{1.3} \int _B e^{i < a, x> } d\mu _{B, h}
(x) = e^{-{h\over 2} |a|^2 },\quad a\in B',\ee
where $< a,
x>$ is the dual product in $B',B$, and where  $|a|$ is the
$H$ norm.
\end{theo}
One may find in standard references the exact  properties on the space $B$; a triplet
$(i,H,B)$ is said a Wiener space. The $L^2$ space of our Weyl
operators will be $L^2 (B, \mu _{B , h/2})$.
Since $h$ stands for the semiclassical parameter, the variance is then $h/2$.
This space is isomorphic to the Fock space ${\cal
F}_s ( H_{\bf C})$ (see  \cite{Sim-PP2}, Theorem I.11,
\cite{Jan}, Theorem 4.1) through the Segal isomorphism.

As in \cite{AJN-JFA} we start by associating with $F\in S(H^2 , Q_A)$ a
quadratic form on a dense subspace of $L^2 (B, \mu _{B , h/2})$,
namely the stochastic extensions of polynomial functions on $H$ (see Section
\ref{s2}). This quadratic form is named $Q_h^{weyl} (F)$ and it
depends on $h>0$. Such a process extends to functions which are not
in $S(H^2 , Q_A)$ such as continuous linear forms on the phase space
$H^2$; in finite dimension $d$, on weaker conditions one associates with
symbols, operators from ${\cal S} (\R^d)$ into ${\cal S}' (\R^d)$, but
in order not to enter in duality we had rather speak of quadratic
forms.

Thereafter in Section \ref{s3}, we prove that the quadratic form
$Q_h^{weyl} (F)$ is related to a bounded operator in $L^2 (B, \mu _{B ,
h/2})$ whose norm is estimated in term of $Q_A$.
We need Definition \ref{d1.3}.

\begin{defi}\label{d1.3}
Let $E$ be an Hilbert space, $q$ a scalar product on $E$, $T\in{\cal
L}(E)$, $T^*=-T$ for $q$, being trace class for $q$. Set $|T|_q$
the unique operator in ${\cal L}(E)$, selfadjoint for $q$ and positive, satisfying $|T|_q ^2 = - T^2$. For $z\in {\bf C}$, we set, \be\label{1.4} D(
E , q, z , T) = {\rm det} (I + z |T|_q ),\ee
where ${\rm det}$ is the Fredholm determinant. Finally we denote by
$\Vert T \Vert _{E , q }$ the $q(X , X)^{1/2}$ norm of $T$.
\end{defi}

In this work, the configuration space will be a real separable Hilbert space $H$.
The natural norm in $H$ or in $H^2$ will be denoted by $|\cdot |$ and the scalar product
of $a$ and $b$ by $a\cdot b$. We consider $A\in {\cal L}(H^2)$, $A=A^*$, nonnegative
trace class.  The QED applications lead us to assume that
$Q_A$ is only nonnegative. Then we define in $H^2$ a quadratic form $Q_A$ by (\ref{1.2})
and a second scalar product by $q_A (X , Y) = (AX) \cdot Y$.
Now we shall define a Hilbert space $B_A$, with a scalar product and  an operator
as in Definition \ref{d1.3}  for the expression of the norm in the first main theorem.

 \begin{prop}\label{p1.5}
Let $H$ be an  Hilbert space, $A\in{\cal L}(H^2)$, selfadjoint, nonnegative and trace class in $H^2$, and let
 $q_A(X , Y ) = (AX)\cdot Y$ and $Q_A(X) = (AX)\cdot X$.  Let $B_A$ be the completed Hilbert space of
 $H^2/{\rm Ker} A$ for $Q_A ^{1/2}$ and
 ${\cal F}: (x , \xi)\mapsto (-\xi, x)$. Then ${\cal F} A$,
 defined in $H^2/{\rm Ker} A$,
extends to a trace class operator in  $B_A$, skew symmetric  for $q_A$; this extension is still denoted  by ${\cal F} A$.
\end{prop}

  This proposition will be proved in Section \ref{s3.2}. The operator ${\cal F} A$ is the fundamental matrix. In H\"ormander \cite{Hor}, the trace $ |{\cal F} A | _{q_A} $
is denoted by $2 {\rm Tr}^+ (A)$.
Moreover, in the finite dimensional case, the dual form of $Q_A$ with respect to the symplectic form is,
$$Q_A^{\sigma} (X) = \sup_Y {| ( ({\cal F}X)\cdot Y |^2 \over
 Q_A(Y)}
  = ( A^{-1}{\cal F}  X ) \cdot ({\cal F}  X ),  $$
and so the uncertainty principle $Q_A \leq Q_A^{\sigma}$ of H\"ormander reads as $\Vert {\cal F} A \Vert _{H^2 ,
q_A } \leq 1$.

We are now able to state,

\begin{theo}\label{t1.6} Let $H$ be a real separable Hilbert space,  $A$ a symmetric nonnegative trace class operator  in $H^2$. Let $B_A$ be the completion of
$H^2/{\rm Ker} A$ for  $Q_A ^{1/2}$. For $F\in S(H^2 , Q_A)$, let $Q_h ^{weyl}(F) $ be the quadratic form of Definition \ref{d2.6}. Then, there exists a bounded operator in
$L^2(B , \mu _{B , h/2})$ denoted by $Op_h ^{weyl} (F) $ satisfying for all $f$ and $g$ in the domain of $Q_h
^{weyl}(F) $,
\bl{1.6}Q_h ^{weyl } (F) (f , g) = < Op_h ^{weyl} (F )f , g>.
\el
Moreover, with the notations of Definition \ref{d1.3} and Proposition \ref{p1.5},
\bl{1.7}\Vert  Op_h ^{weyl} (F ) \Vert  \leq \Vert F   \Vert _{  Q_A }
 D(B_A ,  q_A , 81 \pi  h \max (1, \Vert {\cal F} A \Vert _{B_A , q_A
 } ) \ ,   \ ({\cal F} A) )^{1/2}.
\el
\end{theo}

When the dimension of $H$ is finite, and denoted by $d$, we need only a finite number
of derivatives. In this case, $B=H$ and $L^2 (H, \mu _{H, h/2} )$ is isomorphic to
 $L^2 (H, \lambda )$, where $\lambda $ is the Lebesgue measure, and our operators
 are the usual semiclassical PDO. In this case,
the estimate (\ref{1.7}) may be replaced by,
\be\label{1.5a} \Vert Op_h^{weyl} (F ) \Vert \leq \Vert F \Vert
_{4d, Q_A }
 D \Big (H^2/{\rm Ker} A\ ,  \ q_A\  ,\  81 \pi  h  \max (1, \Vert
 {\cal F} A \Vert _{H^2 , q_A } ) \ ,
 \ ({\cal F} A)
 \Big ) ^{1/2}.
\ee

We now turn to a Beals type characterization result. Beals type theorem consists in proving the existence of an isomorphism between families of functions $(F_h)$ and families of bounded operators $(A_h)$ via the maps
$F_h \mapsto Op_h^{weyl} (F_h) $.

For $V = (a , b)\in H^2$,
we denote by $L_h V$ the operator $L_h V = Op_h ^{weyl} (F)$ with
$F(x, \xi) = - b \cdot x + a \cdot \xi$ (it is proved in $(\ref{2.6})$ below that $L_h V$ has a meaning). For $A$
bounded in  $L^2(B , \mu _{B , h/2})$ and $V_1$, \dots  ,$V_m$
in $H^2$, the iterated commutator  $ {\rm ad} (
L_h V_1) \dots  {\rm ad} ( L_h V_m) A$ is well defined as a quadratic form on the space of stochastic extensions of polynomial functions. When $A= A_h =  Op_h ^{weyl} (F)$ with $F\in S(H^2, Q)$,
$$ {\rm ad} ( L_h V_1) \dots  {\rm ad} ( L_h V_m) Op_h ^{weyl} (F) =
 (h/i)^ {m} Op _h^{weyl } ( d^mF(\,\cdot\, ) (V_1 , \dots  V_m )) $$
 and this commutator is in view of Theorem \ref{t1.6} a bounded operator in
 $L^2(B , \mu _{B , h/2})$ satisfying,
 $$ \Vert {\rm ad} ( L_h V_1) \dots  {\rm ad} ( L_h V_m) Op_h ^{weyl}
 (F)  \Vert
 $$
 $$\leq h^m \Vert F \Vert _{Q_A} \
   D(B_A ,  q_A , 81 \pi  h \max (1, \Vert {\cal F} A \Vert _{B_A ,
   q_A } ) \ ,
    \ ({\cal F} A) )^{1/2}
  \prod _{j= 1}^m Q_A( V_j )
  $$
 and when ${\rm dim\,H}=d<+\infty$ one may replace in the above right hand side
 $\Vert \cdot \Vert _{Q_A} $ by $\Vert \cdot \Vert _{m+ 4d,
 Q_A} $.  So we introduce,

 \begin{defi}\label{d1.7} Let  ${\cal L} (Q) $ be
 the set of operators $(A_h) $
 bounded in $L^2 ( B , \mu _{B , h/2})$  such that  there exists
 $C>0$ satisfying for all $V_1,
 \dots , V_m$ in
 $H^2$,
\bl{1.8}\Vert {\rm ad} ( L_h V_1) \dots  {\rm ad} ( L_h V_m) A_h  \Vert
 \leq
 C \ h^m  \prod _{j= 1}^m Q( V_j ) ^{1/2}.
 \el
 The best above constant $C$  will be denoted by
 $\Vert A \Vert _{{\cal L} (Q)} $. Moreover, for any integer $p$, let
 ${\cal L} _p(Q) $  such that  (\ref{1.8}) holds true for any  $m\leq p$ and let $\Vert A \Vert _{ {\cal L} _p(Q)}$ be the best constant such that  (\ref{1.8}) holds  for any  $m\leq p$.
\end{defi}

With this convention,
\bl{1.9} \Vert  Op_h ^{weyl} (F_h)  \Vert _{{\cal L} (Q_A)}   \leq
 \Vert F \Vert _{Q_A} \
 D(B_A ,  q_A , 81 \pi  h \max (1, \Vert {\cal F} A \Vert _{B_A , q_A
 } ) \ ,
    \ ({\cal F} A) )^{1/2}\el
and when ${\dim\,H}=d<+\infty$, for any integer $m$,
 \bl{1.10}    \Vert  Op_h ^{weyl} (F_h)  \Vert _{{\cal L}_m (Q_A)}   \leq
 \Vert F \Vert _{m+4d, Q_A}  \
 D \Big (H^2/{\rm Ker} A\ ,  \ q_A\  ,\  81 \pi  h  \max (1, \Vert
 {\cal F} A \Vert _{H^2 , q_A } ) \ ,
 \ ({\cal F} A)
 \Big )^{1/2}.\ee

Beals Theorem idea is to prove a converse property with an estimate of the norm.

\begin{theo}\label{t1.9} Set  a real Hilbert space $H$. Let  $A\in{\cal L}(H^2)$ be symmetric and nonnegative and set $Q_A$ the quadratic form (\ref{1.2}). Let $(A_h)$ be a family in  ${\cal L}
 (Q_A) $. Then there exists  $(F_h)$ in $S (H^2 , Q_A)$  such that
 $A_h = Op_h ^{weyl }(F_h)$ and satisfying,
\bl{1.12} \Vert F_h \Vert _{  Q_A} \leq \Vert A_h \Vert _{ {\cal L} (Q_A) }
 \   D(B_A ,  q_A , K  h \max (1, \Vert {\cal F} A \Vert _{B_A , q_A
 } ) \ ,   \ ({\cal F} A) )^{1/2},
 \el
for some $K>0$  universal constant.
\end{theo}

When the dimension of $H$ is finite, and denoted by $d$, the estimation (\ref {1.12}) may
be replaced by:

\bl{1.11a} \Vert F_h \Vert _{ m, Q_A} \leq \Vert A_h \Vert _{ {\cal L}
 _{m+4d} (Q_A) }
 \   D(H^2/{\rm Ker} A\ ,  \ q_A\  ,\   K  h \max (1, \Vert {\cal F}
 A \Vert _{H^2 , q_A } ) \ ,
 \ ({\cal F} A) )^{1/2},\el

See Section \ref{s5.2}.

Now we consider the composition of operators.

\begin{theo}\label{t1.11}
Let $H$ a real Hilbert space, $A\in{\cal L}(H^2)$ symmetric, nonnegative and $Q_A$ the quadratic form (\ref{1.2}) and let $(F_h)$,
 $(G_h)$ two families of symbols in   $S  (H^2 , Q_A) $. Then there exists $(K_h)$ in $S(H^2 , 4Q_A) $ satisfying
\bl{1.13} Op_h ^{weyl} (F_h) \circ  Op_h ^{weyl} (G_h) =  Op_h ^{weyl}
(K_h)\el
    and
\bl{1.15} \Vert K_h \Vert _{ 4Q_A } \leq \Vert F_h   \Vert_{ Q_A  } \
 \Vert G_h   \Vert_{ Q_A  }\
 D(B_A ,  q_A , 81 \pi  h \max (1,  \Vert {\cal F} A \Vert _{B_A ,
 q_A } ) \ ,   \  ({\cal F} A) )
\el
 $$\dots  D(B_A ,  4 q_A , K  h \max (1, 4 \Vert {\cal F} A \Vert _{B_A
 , 4q_A } ) \ ,   \ 4 ({\cal F} A) )
      ^{1/2},$$
where $K$ is a universal constant.
\end{theo}

 When the dimension of $H$ is finite and denoted by $d$, the estimation (\ref{1.15}) may be
 replaced by the following, if $m\geq 8d$,
 \bl{1.14} \Vert K_h \Vert _{ m-8d , 4Q_A} \leq \Vert F_h   \Vert_{m, Q_A } \
 \Vert G_h   \Vert_{m, Q_A  }
 D \Big (H^2/{\rm Ker} A\ ,  \ q_A\  ,\  81 \pi   h  \max (1,  \Vert
 {\cal F} A \Vert _{H^2 , q_A } ) \ ,
 \ ({\cal F} A)
 \Big ) \ee
 $$ \dots   D \Big (H^2/{\rm Ker} A\ ,  \ 4q_A\  ,\  K   h  \max (1, 4
 \Vert {\cal F} A \Vert _{H^2 , 4q_A } ) \ ,
 \ (4{\cal F} A)
 \Big )^{1/2}, $$
It is possible, although not done here, to write a full asymptotic expansion (in term of power of $h$) for the symbol of the composition of two observables.

The interaction between the quantized electromagnetic field and a finite system of fixed  spin $ 1/2$ particles  will use Weyl operators related to the previous class $S(H^2, Q)$; the configuration space $H$ will be
 $H=\{ f\in L^2 (\R^3, \R^3)\ |\ k\cdot f(k) = 0,\ k\ {\rm a.e.}\}$. Photons in vacuum are modelized by an Hamiltonian $H_{ph}$ in a   Hilbert space ${\cal
 H}_{ph}$; the particles system is modelized by a finite dimensional Hilbert space ${\cal H}_{sp}$.
 The full Hamiltonian is a selfadjoint operator $H(h)$ in  ${\cal H}_{ph}\otimes {\cal H}_{sp}$ (see Section \ref{s7.2}).

 Let  $A$  be an observable, i.e., a selfadjoint operator bounded or not in $H_{ph}\otimes H_{sp}$. One may define
 \bl{1.16}  A_h (t) = e^{  i {t\over h} H(h)}    A_h
 e^{ - i {t\over h} H(h)},\quad
 A_h^{free}  (t) =  e^{  i {t\over h} (H_{ph} \otimes I) }    A_h
 e^{ - i {t\over h} (H_{ph} \otimes I)}.\el
 The operator $e^{  i {t\over h} (H_{ph} \otimes I) } $ is a metaplectic operator related to some symplectic map in $H^2$; this symplectic map
 $\chi_t$ belongs to the group studied by B. Lascar \cite{Las}.

  The general framework of Section \ref{s7} is that the full evolution of the usual observables may be described as a perturbation of the free evolution and the perturbation at time $t$ is defined with an operator with a symbol belonging to  $S(H^2, Q_t)$ where $Q_t$  is a nonnegative time dependent quadratic
form defined on the phase space $H^2$ in (\ref{7.14}).

For instance, if $A_h$ is one of the $B_j (x) $ ($1\leq j \leq
 3$) being the component of the quantized magnetic field at  a point $x\in \R^3$, we prove that (Theorem \ref{t7.2}),
$$ B_j  (x , t,   h)=  B_j ^{free} (x , t, h) +  h Op_h ^{weyl} ( B_j
^{res} (x , t, \cdot , h) ), $$
where $(q , p) \mapsto B_j ^{res} (x , t, q, p ,
h)$ is a map from  $S(H^2 , 4 Q_t)$ taking values in ${\cal L} ({\cal
H}_{sp})$. The same holds true for the electric field and the spin observables. The photon number evolution is defined in Theorem \ref{t7.3}.

\section{The Weyl calculus.}\label{s2}

Let us start with the finite dimension case. To any function $F$, say $F\in S(H^2, Q)$, we associate an operator $Op_h^{weyl}(F)$. We have, for any functions $f$ and $g$ defined on $H$, for example $C^\infty$ with compact support, using the Gaussian measure instead of the Lebesgue measure in order to prepare the infinite dimension case,
\bl{2.1} < Op_h^{weyl}(F) f , g> = \int _{H^2} F(Z) H_h^{gauss} (f , g) (Z)
d\mu _{H^2 , h/2} (Z),
\el
 where $ H_h^{gauss} (f , g)$ is the Wigner function which is  directly  defined below in infinite dimension.

\subsection{Wigner function.}\label{s2.1}

We first introduce Weyl translations associated to $X$ in the phase space $H^2$. In finite dimension, these Weyl translations are defined for example in \cite{CR}.

For any $a \in B' \subset H$, the function $B\ni x \rightarrow a(x)
$ belongs to $L^2(B , \mu _{B , h/2})$ with a norm  being $ \sqrt { h/2} |a|$ (where $|a|$ is the norm of $a$ in $H$), so the map which associates to $a\in B'$ this $L^2$ function has by density an extension from $H$ into $L^2(B , \mu
_{B , h/2})$ denoted by $a \mapsto \ell _a$.

For each $X=(a , b)\in H^2$, the operator $V_h(X)$ (Weyl translation of $X$) is unitary in $L^2(B , \mu _{B , h/2})$, it is defined by, 
\bl{2.2} ( V_h(X) f ) (u) = e^{{1\over h} \ell _{a+ib} (u) -{1\over 2h}
|a|^2
- {i\over 2h} a\cdot b } f (u - a),\quad f\in L^2(B , \mu _{B ,
h/2})\el 
and it may be viewed as an operator in the Fock space ${\cal F}_s (H_{\bf C})$ through the Segal isomorphism.

\begin{defi}\label{d2.1}  For $f$ and $g$ in $L^2 (B ,
\mu _{B , h/2})$ or in ${\cal F}_s (H_{\bf C})$, we define the Wigner function $ H_h^{gauss}(f , g ) $  on $H^2$ by, 
\bl{2.3} H_h^{gauss}(f , g ) (Z) =  e^{|Z|^2 \over h}   < V_h(-2Z) f ,
 \check g >,\quad  Z\in H^2,\el
where $\check g (u) = g(-u)$.
\end{defi}

Since $V_h(-2Z)$ is unitary, one has,
\bl{2.4} |H_h^{gauss}(f , g ) (Z)| \leq e^{{1\over h}|Z|^2 }
 \Vert f \Vert _{ L^2 (B , \mu _{B , h/2})} \
 \Vert g \Vert _{ L^2 (B , \mu _{B , h/2})}.\el
The Wigner function is also equivalently defined by, 
 $$ H_h^{gauss}(f , g ) (z, \zeta) = e^{|\zeta |^2 \over h} \int _B
 e^{-2 {i\over h} \ell _{\zeta } ( t ) } f (z +t) \overline { g
 (z-t)}
 d \mu _{B , h/2} (t),\quad Z = (z, \zeta) \in H^2.  $$

\subsection{Stochastic extensions.}\label{s2.2}

When ${\rm dim}\,H=\infty$, the integral (\ref{2.1}) is meaningless since $F$ and  $H_h^{gauss}(f , g )$ are defined only on $H^2$ whereas one has to integrate (for the Gaussian measure) on $B^2$.

However, for some $F$, it is possible to construct a convenient extension  $\widetilde F$ on
$B^2$. This has been achieved by L. Gross in \cite{Gro} who gave to $\widetilde F$ the name of stochastic extension. This construction is recalled in  \cite{AJN-JFA} in the context of extension in $L^p$ (whereas it is rather a convergence in measure in \cite{Gro}). We remind here the main features.

If $E$ is a $n$ dimensional subspace of $H$, we set
 $\widetilde \pi _E : B \rightarrow E$ the map defined by,
\bl{2.5} \widetilde \pi _E  (x) = \sum _{j=1}^n \ell _{u_j}(x) u_j
\quad{\rm a.e.\ in}\ B,\el
where $u_j$, $1\leq j \leq n$, is an orthonormal basis of $E$.
The map $\widetilde \pi _E $ is independent of the choice of the basis.

\begin{defi}\label{d2.2} Let $(i, H, B)$ be a
Wiener space,   set $p\in [1 , \infty )$ and $h>0$. We say that a Borel function $F$ on  $H^2$ has a stochastic extension $\widetilde
F$ in $L^p (B , \mu_{B , h})$ if for any increasing sequence $(E_n)$ of finite dimension subspaces of $H^2$ with dense union, the functions $F \circ \widetilde \pi _{E_n}$ are in $L^p (B , \mu_{B
, h})$ and if $F \circ \widetilde \pi _{E_n}$ converges to
$\widetilde F$ in $L^p (B , \mu_{B , h})$.
\end{defi}

The fact that a function $F$ on an Hilbert space $H$ has or not a stochastic extension in the above sense is independent of the chosen space $B$. In this sense, the choice of $B$ is irrelevant. Indeed, let  $(E_n)$ be an increasing sequence of subspaces of $H$ and set for  $m <n$,  $\pi _{mn}: E_n \rightarrow E_m$ the orthogonal
projection. When  $f$ is a Borel function we set,
$$ I_{mn}(f) = \int _{E_n} | f(x) - f( \pi_{mn} (x) |^p d\mu_{E_n ,
h}(x), $$
and our definition concerning stochastic extensions means that $I_{mn}(f)\rightarrow 0$ when  $ m = \inf (m,
n)\rightarrow\infty$.  For $p=1$, if $(i_1 , H ,
 B_1)$
and $(i_2 , H , B_2)$ are Wiener spaces and if $f$ has a stochastic
extension in $L^1$ for the variance $h$ then the two stochastic extensions of $f$ have in  $B_1$ and $B_2$ the same integrals.

Many classes of functions have stochastic extensions, see examples in \cite{AJN-JFA}. We shall see that the symbol $F$ and the Wigner function has stochastic extensions for different reasons. We begin with the symbol $F$, being in  $S(H^2,
Q_A)$.

\begin{theo}\label{t2.3}(See \cite{Jag}, Proposition 3.10) Set $h>0$ and  $p\in[1,+\infty)$. Let
$A\in {\cal L}(H^2)$ be a nonnegative selfadjoint and trace class operator in  $H^2$.
Then, any $f\in S(H^2, Q_A)$ has a stochastic extension $\widetilde f $ in $L^p(B,\mu_{B,h})$, which is bounded almost everywhere by
 $||f||_Q$.

In addition, if
  $(u_j)$  is an orthonormal basis of $H$ constituted with eigenvectors of $A$  and if the corresponding eigenvalues are the  $\lambda _j $, if $E$ is a finite dimension subspace of $H$, one has,
$$
  ||f \circ \tilde{\pi}_{E} - \widetilde{f} (x) ||_{L^p(B, \mu _{B ,
  h })}
\leq C(p) h^ {1/2}  \Vert f \Vert _{Q_A} \left(\sum_{j\geq 0}\lambda
_j |\pi _{E} (u_j )-u_j|^{\alpha(p)} \right)^{1/\alpha(p)},
$$
where
$$
C(p)=K(p)\left( \sum_{0}^{\infty} \lambda_j \right)^{(1/2)
-(1/p)},\  \alpha(p)= p, \  {\rm for} \  p\geq 2,
$$
$$C(p)= 1,\ \alpha(p)= 2,  \ {\rm for} \ p\leq 2$$
and where $
 K(p)= 2^{1/2} \pi^{-1/2p}
\left( \Gamma \left(   {p+1 \over  2} \right)\right)^{1/p}.$
\end{theo}

Note that if $f,g\in L^2(B , \mu _{B , h/2})$ then their Wigner function has in general no stochastic extension in $L^2 (B^2 , \mu _{B^2 , h/2})$.

However, any polynomial function has a stochastic extension in $L^p (B^2 , \mu _{B^2 ,
h/2})$. This is proved in Proposition 3.17 in \cite{Jag}. This point motivates the following.

\subsection{Weyl quadratic form.}\label{s2.3}

\begin{defi}\label{d2.6} Let $F$ be a Borel function on $H^2$ with a stochastic extension $\widetilde F$  belonging to $L^2 (B^2 , \mu _{B
 ^2, h/2 })$. For  $f$ and $g$ polynomial functions, we set,
\bl{2.6} Q _h^{weyl } (F) (f , g)  = \int _{B^2} \widetilde F (Z)
\widetilde H_h^{gauss} ( f , g) (Z) d\mu _{B^2 , h/2} (Z), \el
where $\widetilde H_h^{gauss} ( f , g)$ is the stochastic extension of the Wigner function  $H_h^{gauss} ( f , g)$ in $L^2 (B^2 , \mu _{B ^2, h/2 })$.
\end{defi}

We remind that the Wigner function of  two polynomial functions is also a polynomial function. Consequently, this Wigner function has also a stochastic extension $\widetilde F$ in $L^2 (B^2 , \mu
 _{B ^2, h/2 })$. More generally, we may define $Q _h^{weyl } (F) (f , g)$  as soon as $(f,g)$ is  such that  $H _h^{gauss} (f , g)$ has a stochastic extension in $L^2 (B^2 , \mu _{B ^2, h/2 })$.

In view of Theorem \ref{t2.3}, if  $F$ belongs to $S(H^2, Q_A)$ then $F$ has  a stochastic extension  $\widetilde F$ in $L^2 (B^2 , \mu
 _{B ^2, h })$ which is bounded. Consequently, if  $f$ and $g$ are two polynomials functions then,
\bl{2.7} |Q _h^{weyl } (F) (f , g) | \leq\Vert F \Vert _{\infty} C(f , g),\quad  C(f , g) = \Vert \widetilde H_h^{gauss} ( f , g) \Vert _{
L^1 (B^2 , \mu _{B^2 , h/2})}. \el
If $F$ is a continuous linear form on $H^2$ then $F$ also has a stochastic extension in  $L^2 (B^2 , \mu _{B ^2, h })$ and $Q _h^{weyl } (F)$ is well defined.

\subsection{Wick symbol.}\label{s2.4}

For any $X =(a , b)\in H^2$, the coherent state $
\Psi_{(a , b) , h}$ is defined as,
\bl{2.8}  \Psi_{X , h} (u) = e^{{1\over h} \ell _{ (a+ib)} (u)  -{1\over
2h}|a|^2 -
{i\over 2h} a\cdot b},\quad X = (a , b) \in H^2,\quad{\rm a.e.}\ u\in B.\el

\begin{prop}\label{p2.7} For $X$ and $Y$ in $H^2$, one has
\bl{2.9}H_h^{gauss} ( \Psi_{X , h}  , \Psi_{X , h} ) (Z)=
e^{- {|X|^2 \over h} + {2\over h} X \cdot Z }.\el \end{prop}

{\it Proof.} \newline
In view of $\Psi_{X, h} = V_h(X) \Psi_{0,
h}$ and $\check \Psi_{X, h} = \Psi_{-X, h} $,
$$ V_h(X) V_h (Y)   = e^{{i\over 2h}\sigma (X , Y)}  V_h (X+Y). $$
Thus, for any  $X$ and $Z$ in $H^2$, one sees,
\bl{2.10} \Psi_{X+ Z, h} = e^{{i\over 2h} \sigma (X , Z) } V_h(Z) \Psi_{X,h}.\el
Besides,
\bl{2.11} < \Psi_{X h} , \Psi _{Yh}> =e^{-{1\over 4h}(|X-Y|^2 ) + {i\over 2h} \sigma (X , Y)}.
\el One then deduces (\ref{2.9}).

\hfill$\Box$

As a consequence, if $X$ and $Y$ are in $H^2$ then the Wigner function $H_h^{gauss} ( \Psi_{X , h}  , \Psi_{X , h} )$ admits a stochastic
extension
\bl{2.12} \widetilde H_h^{gauss} ( \Psi_{X , h}  , \Psi_{X , h} ) (Z)=
e^{- {|X|^2 \over h} + {2\over h} \ell _X ( Z)}.\el
For any $F$ in $S(H^2 , Q_A)$ and each $X,Y$ in
$H^2$,  $Q_h^{weyl} (F)( \Psi_{X , h}  ,
\Psi_{X , h} )$ is defined and so is the Wick symbol of $Q_h^{weyl} (F)$ defined by,
\bl{2.13}
\sigma _h ^{wick} (Q_h^{weyl} (F) ) (X) = Q_h^{weyl} (F)( \Psi_{X ,
h}  , \Psi_{X , h} ),\quad X\in H^2 .\el

The Wick symbol is related to $F$ via the heat operator, as in finite dimension. Set
 $F\in S(H^2, Q_A)$ and $X\in H^2$. Also set,
$\tau _X (F)\, :\, Y \mapsto F(X+Y)$. Then $\tau _X (F)\in S(H^2,
Q_A)$ and let $\widetilde \tau _X F $ be its stochastic extension in
 $L^1 ( B^2 , \mu_{B^2, h/2})$. The heat operator is defined by,
\bl{2.14}(H_{h/2} F) (X) = \int _{B^2} \widetilde \tau _X F (Y)d\mu_{B^2,h/2} (Y).\el
For $E\subset H^2$ with ${\rm dim\, E}<+\infty$, one may define,
\bl{2.15} (H_{E , h/2} F) (X) = \int _{E}   F(X+Y)  d\mu_{E, h/2}
(Y).\el
If $F\in S(H^2,Q_A)$ then  these two functions are in $S(H^2,Q_A)$. It is proved in \cite{AJN-JFA} (Theorem 7.1)  that, if $F\in S(H^2,Q_A)$,
\bl{2.16}\sigma _h ^{wick} ( Q_h^{weyl} (F)) (X) = (H_{h/2} F) (X).
\el
We note that this point is also a consequence of Proposition \ref{p2.7} and of
(\ref{2.12}). Indeed, using translation change of variables for Gaussian measures (see for example Proposition 4.2 in \cite{AJN-JFA}), one has
$$ Q _h^{weyl } (F) (\Psi_{X h} ,\Psi_{X h} )  = \int _{B^2}
\widetilde F (Z)
e^{- {|X|^2 \over h} + {2\over h} \ell _X ( Z) }
 d\mu _{B^2 , h/2} (Z)
 = \int _{B^2} \widetilde \tau _X F (Y)d\mu_{B^2, h/2} (Y)
 = (H_{h/2} F) (X).   $$

\section{Extensions and norm estimates.}\label{s3}

\subsection{The finite dimension case.}\label{s3.1}

We first recall a classical result. A proof is available in  H\"ormander \cite{Hor} (Theorem 21.5.3).

\begin{theo}\label{t3.1}  Set $E$  a real Hilbert space of finite dimension $2d$, with a symplectic $2-$form  $\sigma$. Let   ${\cal
  F}\in {\cal L} (E)$  be such that  $\sigma (X , Y) = <  X , {\cal F} Y>$, for all $X$ and $Y$ in $E$. Let $A\in {\cal L}(E)$ being selfadjoint and nonnegative. Then,
\bl{3.1} (A{\cal F} A X ) \cdot Y = - (A X ) \cdot ({\cal F} A Y), \quad X , Y \in E.\el
Moreover, there exists a symplectic basis $(U_j) , (V_j) $ of $E$ ($j=1,\dots ,d)$ and there exists an integer  $p\leq d$  satisfying,

  i)  We have $(AX) \cdot Y = 0$, if $X\not =Y$, for  $X$ and $Y$ belonging to the basis.

  ii)  We have $AV_j = 0$, for $j\leq p$.

  iii) We have for $j>p$,
\bl{3.2}  < A U_j, U_j > = < A V_j, V_j > = \lambda _j,  \el
 where the $(\lambda _j )^2$ are the non vanishing eigenvalues of  $ - ({\cal F} A)^2$, which are nonnegative. If  $j \leq p$ then one has $({\cal F} A)^2 U_j = ({\cal F} A)^2 V_j = 0$. If $j>p$ then
 $({\cal F} A)^2 U_j = - (\lambda _j )^2 U_j$ and  $({\cal F} A)^2 V_j
 = - (\lambda _j )^2 V_j$.
\end{theo}

\begin{theo}\label{t1.4} Let $H$ be a real Hilbert space of finite dimension
$d$, and $A\in{\cal L}(H^2)$, $A^*=A$, $A$ trace class and nonnegative. Let us set $q_A (X, Y) = (AX) \cdot Y$ and $Q_A (X) = q_A(X
, X)$.
 Let  ${\cal F}$ be defined on $H^2$ by
 ${\cal F} (x , \xi) = (-\xi, x)$.
 Then, for any $F\in S_{4d} ( H^2, Q_A)$,
\be\label{1.5} \Vert Op_h^{weyl} (F ) \Vert \leq \Vert F \Vert
_{4d, Q_A }
 D \Big (H^2/{\rm Ker} A\ ,  \ q_A\  ,\  81 \pi  h  \max (1, \Vert
 {\cal F} A \Vert _{H^2 , q_A } ) \ ,
 \ ({\cal F} A)
 \Big ) ^{1/2}.
\ee
\end{theo}

{\it Proof of Theorem \ref{t1.4}.}
\newline
Let $H$ be a real Hilbert space with finite dimension $d$ and let $E=H^2$, $(e_j)$ an orthonormal  basis of $H$. We set $u_j = (e_j, 0)$ and $v_k = (0, e_k)$. Let  $(U_j , V_k) $ and $p$ be given by Theorem \ref{t3.1}. There exist a symplectic map $\chi$  such that   $\chi (u_j)
= U_j$ and $\chi (v_k)= V_k$. In view of Theorem 18.5.9 in \cite{Hor},
there is a metaplectic transform   $U_{\chi}$, unitary in $L^2 (E, \mu _{E ,h/2})$ satisfying,
\bl{3.3}  U_{\chi}^{\star} Op_h ^{weyl} (F)  U_{\chi} = Op_h ^{weyl}
 (F\circ \chi).\el
For any $F\in S_{4d} ( H^2 , Q_A )$, one has
\bl{3.4} \Vert  Op_h ^{weyl} (F) \Vert  = \Vert U_{\chi} ^{\star}   Op_h
 ^{weyl} (F) U_{\chi} \Vert
  = \Vert  Op_h ^{weyl} (F \circ \chi ) \Vert.\el
The symbol $F \circ \chi$ belongs to $S_{4d}( H^2 , Q_A \circ \chi)$. Let $(x,\xi)$, $x=(x',x'')$, $\xi=(\xi',\xi'')$, $x',\xi'\in\R^p$, $x'',\xi''\in\R^{d-p}$ be the coordinates chosen according to Theorem \ref{t3.1}. Writing  $G=F \circ \chi$, $G$ is a function of $(x',\xi', x'' , \xi'')$. By definition, the fact that  $F \circ \chi$ belongs to
 $S_{4d}( H^2 , Q_A \circ \chi)$ may be written as, with obvious notations, for all multi-indices
 $( \alpha ', \alpha '' , \beta', \beta '' ) $ with a length
 $\leq 4d$,
 $$ | \partial _{x'} ^{\alpha '} \partial _{x''} ^{\alpha ''}
 \partial _{\xi'} ^{\beta '}
 \partial _{\xi''} ^{\beta ''}
 G(x',  \xi',x'',  \xi'')|
 \leq \Vert F \circ \chi  \Vert _{ 4d, Q_A \circ \chi }\
 \prod _{j}  Q_A ( U'_j) ^{(1/2) \alpha '_j} Q_A ( U''_j) ^{(1/2) \alpha ''_j}
 Q_A ( V'_j) ^{(1/2) \beta '_j}
 Q_A ( V''_j) ^{(1/2) \beta ''_j}. $$
Since $AV_j = 0$ for $j\leq p$, the above right hand side is vanishing if  one of the $\beta '_j \not  = 0$.
Consequently, $G$ is independent of $\xi '$ and may be denoted by
 $G( x' , x'', \xi '')$.  Since the  $U_j $ and $ V_j$, $j > p$  satisfy (\ref{3.2}), then this function satisfies, when $|\alpha '' | + |\beta '' | \leq 4d$,
$$|  \partial _{x''} ^{\alpha ''} \partial _{\xi''} ^{\beta ''}
 G(  x', x'',  \xi'')|
 \leq \Vert F   \Vert _{ 4d, Q_A  }\  \prod _{j} \lambda _j ^{ (1/2) (
 \alpha ''_j
 + \beta '' _j ) },$$
where the $(\lambda _j )^2$ are the non vanishing eigenvalues of   $ - ({\cal F} A)^2$.
In particular this estimate is true for any
$(\alpha '', \beta '')$ with $\alpha'' _j \leq 2$ and $\beta ''_j \leq 2$ for all $j$. For any $x'$, according to Theorem 1.4 in \cite{AJN-JFA}
applied to the map  $(x'' , \xi'') \mapsto G(x' , x'' , \xi'')$ with the sequence $\varepsilon_j = \sqrt {\lambda _j}$, one has
 $$ \Vert  Op_h ^{weyl} (F \circ \chi ) \Vert  \leq \Vert F   \Vert
 _{4d,  Q_A }
    \prod _{j=p+1}^d   ( 1 + 81 \pi h S \lambda _j ), $$
with $ S = \sup_j {\rm max } (1 , \lambda _j)  $.
One then deduces that,
 $$ \Vert  Op_h ^{weyl} (F ) \Vert  \leq \Vert F   \Vert _{4d,  Q_A
 }
  \prod _{j=p+1}^d   ( 1 + 81 \pi h S \lambda _j ).    $$
Since the $\lambda _j$ are the non vanishing eigenvalues of the operator $|{\cal F}A|_{q_A}$ (see Definition \ref{d1.3}) each computed twice then,
 $$
  D \Big (H^2/{\rm Ker} A\ ,  \ q_A\  ,\  81 \pi  h  \max (1, \Vert
  {\cal F} A \Vert _{H^2 , q_A } ) \ ,
 \ ({\cal F} A)
 \Big )^{1/2}= \prod _{j=p+1}^d  ( 1 + 81 \pi h S \lambda _j ).$$
 \hfill$\Box$

\subsection{The infinite dimension case.}\label{s3.2}

The above Proposition \ref{p1.5} is a consequence of the next result. When $B$ is an Hilbert space, let us denote by $ \Vert \cdot  \Vert _{ {\cal L}^2(B)}$ (resp. $
\Vert \cdot  \Vert _{ {\cal L}^1(B)}$ ) the Hilbert Schmidt norm
(resp. the trace class norm ) of bounded operators in $B$.

 \begin{prop}\label{p3.2} Let $A$ be a selfadjoint nonnegative trace class operator in $H^2$ and define $B_A$ as the completed space of
  $H^2/{\rm ker}A $ with respect to the $Q_A^{1/2}$. Then

 i) For any $T\in {\cal L} (H^2)$ ,
 $T  A ^{1/2}$ is
 defined in $H^2 /{\rm ker }A$ and has an Hilbert Schmidt extension in $B_A$ satisfying,
 \bl{3.5} \Vert T A^{1/2} \Vert _{ {\cal L}^2(B_A)}^2 =
   {\rm Tr}_{H^2} ( A TT^{\star}  ).\el

 ii) Let $E$ be a finite dimension subspace of $H$. Define the orthogonal projection $P_E : H^2/{\rm ker}A \rightarrow E^2/{\rm ker}A $  (with respect to
 $q_A$). Also define the orthogonal projection $\pi_E : H^2 \rightarrow E^2$
 (with respect to
 the usual scalar product of $H^2$). Then, the operator  $\pi_E {\cal F} A P_E $ maps  $
 H^2 /{\rm ker}A $ into itself and it is a skew symmetric operator (for the scalar product $q_A$).

 iii) The two operators ${\cal F} A $ and $\pi_E {\cal F} A P_E $  have trace class skew symmetric extensions  in  $B_A$, for the scalar product $q_A$.

 iv)  Let $(E_n)$  be an increasing sequence of finite dimensional subspaces of  $H$ with  a dense union in $H$. Then,
\bl{3.6} \lim _{n \rightarrow \infty } \Vert {\cal F}  \pi_{E_n}  A
  P_{E_n} - {\cal F} A \Vert
  _{ {\cal L} ^1 ( B_A , q_A) } = 0.\el
 ($ \Vert \cdot \Vert  _{ {\cal L} ^1 ( B_A , q_A) }$ being the trace class operators in  $B_A$ norm).
\end{prop}

{\it Proof. }\newline
 i) Let
  $(u_j)$ an Hilbert basis of  $H^2/{\rm ker}A $ being eigenvectors of $A$ and $\lambda _j $ the corresponding eigenvalues. Then ($v_j) = (u_j / \sqrt {\lambda _j})$
  is an Hilbert basis of $B_A$.
For any operator $T\in {\cal L} (H^2)$, the sequence
  $\Vert T A ^{1/2} v_j \Vert _{B_A}^2 = < T^{\star} A  T u_j ,
   u_j >$ is summable since  $A$ is of trace class. Thus,
  $T  A ^{1/2}$ well defined on $H^2 /{\rm ker }A$ extends to an  Hilbert Schmidt operator in $B_A$ and
  $$ \Vert T A^{1/2} \Vert _{ {\cal L}^2(B_A)}^2 =
  \sum _j \Vert T A^{1/2} v_j \Vert _{B_A}^2 = \sum _j \Vert T  u_j
  \Vert _{B_A}^2
  = \sum _j \Vert  A^{1/2} T  u_j \Vert ^2 = {\rm Tr} (T^{\star} A T
  )
  =  {\rm Tr}_{H^2} ( A TT^{\star}  ). $$

ii) For any $X$ and $Y$  in $H^2$,
$$ (A \pi_E {\cal F} A P_E X )\cdot Y =
- (A P_E X )\cdot ( {\cal F} \pi_E A Y) = - (A  X )\cdot ( P_E {\cal
F} \pi_E A Y) = -  (A  X )\cdot ( {\cal F} \pi_E A Y) $$ 
$$ =   (P_E  \pi_E {\cal F}  A  X )\cdot (   A Y)=
(  \pi_E {\cal F}  A  X )\cdot (   A P_E Y)=
 - (  A  X )\cdot (  \pi_E  {\cal F} A P_E Y).$$
(One uses above that  $P_E$ is selfadjoint for $q_A$,  $A$ and $\pi _E$  are selfadjoint for $\cdot$, ${\cal F} $ skew symmetric for $\cdot$,
 ${\cal F} $ commutes with $P_E$ and $\pi_E$, and $P_E \pi _E = \pi _E$).

iii) The fact that ${\cal F}A$ is skew symmetric for
 $q_A$ is standard. Applying i) with  $T = {\cal F}$ and $T= I$ one proves that
  ${\cal F}  A ^{1/2}$ and $ A ^{1/2}$ have Hilbert Schmidt extensions in  $B_A$. The same argument shows  that $\pi_E {\cal F} A $
has a trace class extension in $B_A$. Besides, $P_E$  extends to a bounded operator in $B_A$ with unit norm. Thus iii) is proved.

iv) One applies i) with  $T =
{\cal F} \pi_{E_n}  -{\cal F} $. According to Theorem 6.3  of Gohberg-Krein \cite{GK} applied in $H^2$,
  $$\lim _{n \rightarrow \infty } {\rm Tr}_{H^2}  \left (  A
  [{\cal F} \pi_{E_n}  -{\cal F}] [{\cal F} \pi_{E_n} -{\cal
  F}]^{\star} \right ) = 0.$$
Thus,
\bl{3.7} \lim _{n \rightarrow \infty } \Vert {\cal F} \pi_{E_n} A^{1/2} -
 {\cal F} A^{1/2}
  \Vert _{ {\cal L}^2(B_A)} = 0.\el
Applying again
Theorem 6.3  of  \cite{GK} but with
 $B_A$, one also has
\bl{3.8} \lim _{n \rightarrow \infty } \Vert  A^{1/2}P_{E_n} - A^{1/2}
 \Vert _{ {\cal L}^2(B_A)} = 0.
 \el
Combining (\ref{3.7}) and (\ref{3.8}), one indeed obtains (\ref{3.6}).

\hfill$\Box$

The proof of Theorem \ref{t1.6} uses the four Propositions \ref{p3.3}-\ref{p3.6} below.

\begin{prop}\label{p3.3} (Ramer) For any finite dimensional $E\subset B'\subset H$,  let $E^{\perp}$ be the orthogonal complement of $E$ in
$H$ and set,
\bl{3.9} \widetilde E^{\perp} = \{ x\in B ,\ \forall\, u \in E, \ u(x) = 0\}.\el
Then $(i, E^{\perp} , \widetilde E^{\perp})$ is a
Wiener space and the corresponding Gaussian measure $\mu _{\widetilde E^{\perp} , t }$ with variance $t>0$ may be identified to the  tensor product measure $\mu _{E, t} \otimes \mu _{\widetilde
E^{\perp} , t }$.
\end{prop}

For  each  function $F\in S(H^2,Q)$ and for any finite dimensional $E\subset B'^2$, we set (with the notations of Proposition \ref{p3.3}),
\bl{3.10} (H_{  E^{\perp}, t } F) (X) = \int _{ \widetilde E^{\perp}}
\widetilde F(X + Y ) d\mu _{\widetilde E^{\perp} , t } (Y),\quad
X\in H^2,\el
 where $\widetilde F$  is a stochastic extension of $F$ (its existence is given by Theorem \ref{t2.3}). The function $H_{  E^{\perp}, t } F$ belongs to $S(H^2, Q)$.

The following Proposition recalls known results in finite dimension but it takes a different form when using Gaussian measures. Replacing $F$ by $H_{h/2}F$ in the formula defining the Weyl calculus allows to take the anti-Wick operator with symbol $F$. Let us recall the definition of the anti-Wick operator. If $H$ is a finite dimension space, if $f$ is a function in $L^2(H,\mu_{H,h/2})$, we denote here by
$H^2\ni X \mapsto T_{X,h}f$ the Segal Bargmann transform of $f$ being defined on $H^2$ by,
$$ (T_{X,h}f)=e^{|X|^2/4h}\int_{H} f(t)\overline{\Psi_{X,h}(t)}d\mu_{H,h/2}(t),$$
where $\Psi_{X,h}$ denotes the coherent state defined in Section \ref{s2}. Then, one associates an anti-Wick operator with any  bounded function $\Phi$ on $H^2$ by,
$$
<Op_h^{AW}(\Phi)f,g>=\int_{H^2} \Phi(X)T_{X,h}f\overline{T_{X,h}g}d\mu_{H^2,h}(X).
$$
It is known that,
$$
\int_{H^2} |T_{X,h}f|^2d\mu_{H^2,h}(X)=||f||^2.
$$
It is also known that the Weyl symbol of $Op_h^{AW}(\Phi)$ is $H_{h/2}\Phi$. If $S$ is a finite dimensional subspace of $H$, replacing $F$ by $H_{S,h/2}F$, amounts to define an operator, called hybrid operator in \cite{AJN-JFA}, being similar to an anti-Wick operator in the variables $S$.
Let $H$ be a finite dimensional space ($d={\rm dim}\,H$), for $X\in H^2$ the Weyl symbol of the orthogonal projection on $\Psi_{X h}$ is
\bl{3.11} F_X (Z) = 2^d e^{-{1\over h} |X-Z|^2},\el
which is meaningless in infinite dimension, contrarily to the Wick symbol  which is  $e^{-{1\over 2h} |X-Z|^2}$.

\begin{prop}\label{p3.4} Let $H$ be a finite dimensional Hilbert space written as $H = E \oplus S$ where $E$ and $S$ are two orthogonal subspaces. Denote by $(X', X'')$ the coordinates in $H^2$ according to the decomposition $H^2 = E^2\oplus S^2$. For each $X''\in S^2$, let $\Psi _{X'', h}$ denotes the corresponding coherent state defined in  Section \ref{s2} replacing $H$ by $S$. This coherent state belongs to $L^2(S , \mu _{S , h/2} )$. For any
polynomial function $f$ and each $X''$ in $S^2$, set
\bl{3.12}  (T_{X'', h} f ) (t') = e^{|X''|^2 /4h} \int _S f(t', t'')
\overline {
\Psi _{X'', h} (t'') } d\mu_{S , h/2} (t''). \el
For any $G \in S( E^2 , Q)$,  $Q_h ^{weyl, E}(G)$  denotes the Weyl quadratic associated with $G$ and  with
$H$ replaced by $E$. Then, for all $\Phi \in S(H^2 ,
Q)$ and for any  polynomial functions $f,g$, we have the following identity,
\bl{3.13} Q_h ^{weyl} (H_{S , h/2}\Phi )  (f, g) = \int _{S^2} Q_h ^{weyl,
E} (\Phi ( \cdot , X'') )
   ( T_{  X'', h }f ,  T_{  X'', h }g) d\mu _{ S^2, h} (X'').\el
\end{prop}

{\it Proof. } \newline
 According to  Proposition 4.5 of \cite{AJN-JFA} one has for any functions $\Phi$ and $G$ on $H^2$,
$$\int _{H^2} (H_{S , h/2} \Phi  ) (X) G(X)d\mu _{ H^2, h/2} (X)=
\int _{H^2}
 \Phi (X', X'') ( M_{S , h/2} G) (X', X'') d\mu _ {E^2 , h/2 } (X')
 d\mu _ {S^2 , h } (X''),$$
where we set,
$$ ( M_{S , h/2} G) (X', X'') = \int _{S^2} G \left ( X', { X'' \over
2 } + Y''
\right ) d\mu _ {S^2 , h/4 } (Y''). $$
One applies this identity with $G(X)= H _h^{gauss} (f , g ) ( X)$. One obtains,
 $$ Q_h ^{weyl} (H_{S , h/2}\Phi )  (f, g) = \int _{H^2}
  \Phi (X', X'') ( M_{S , h/2} H _h^{gauss} (f , g )) (X', X'') d\mu
  _ {E^2 , h/2 } (X')
 d\mu _ {S^2 , h } (X'').$$
If the subspace $S$ is finite dimensional,
$$ ( M_{S , h/2} H _h^{gauss} (f , g )) (X', X'') =  e^{ {1\over 2 h}
|X''|^2 }  \int _{S^2}
 H _h^{gauss} (f , g ) ( X', Y'') F _{ X''} (Y'')
 d\mu _ {S^2 , h/2 } (Y''), $$
with
$  F _{ X''} (Y'') = 2^d e^{ -{1\over h} |X'' - Y''|^2 } $.
Since  $  F _{ X''} (Y'')$ is the Weyl symbol of the orthogonal projection in $L^2 (S , \mu _{S , h/2})$ on $\C\Psi_{X'' , h}$, one deduces that,
$$ ( M_{S , h/2} H _h^{gauss} (f , g )) (X', X'') =  H_h^{gauss , E}
( T_{X'' , h} f,
T_{X'' , h} g ) (X'), $$ where $H_h^{gauss , E}$ is the Wigner function for two functions defined on $E$. The proof of the Proposition \ref{p3.4} then follows.

\hfill$\Box$

\begin{prop}\label{p3.5} Let $S$ be a subspace of $H^2$ being either finite dimensional or which is the orthogonal of a  finite dimensional subspace included in $B'^2$.  Let $\pi : H^2 \rightarrow S$ denotes the corresponding orthogonal projection. Then, for   $F\in S(H^2,Q_A)$ one has,
\bl{3.14}  \Vert H_{ S  , t }F  - F \Vert _{Q_A}  \leq \sqrt {t} \Vert F
\Vert _{Q_A}
 {\rm Tr} (A \pi _S  )^{1/2}.\el
 \end{prop}

{\it Proof.}\newline
 Let $(u_j)$ be an Hilbert basis of $H^2$
 constituted with eigenvectors of $A$ and let $\lambda _j$ be the corresponding eigenvalues.  Using Definition (\ref{2.15}) when $S$ is finite dimensional
or (\ref{3.10})  when the orthogonal complement of $S$ is included in $B'$. In both cases, one has
$$  H_{ S  , t }F (X)  - F (X) = \int _{ \widetilde S}
(\widetilde F(X + Y )- F(X))  d\mu _{\widetilde S , t } (Y),  $$ 
where $\widetilde S = S$ if $S$ is finite dimensional and  where $\widetilde
S $ is defined in (\ref{3.3}) in the other case.

Since $F$ belongs to $S(H^2,Q_A)$ and possibly using its stochastic extension, one has for any  $X\in H^2$ and $Y\in\widetilde S $,
$$ |\widetilde F(X + Y )- F(X)| \leq \Vert F \Vert _Q \widetilde
Q(Y)^{1/2}, $$
where $\widetilde Q = Q$ if $S$ is finite dimensional and in the other case,  $\widetilde Q$ is the stochastic extension of $Q$  which exists from \cite{AJN-JFA} and satisfies
$$ \widetilde Q(Y) = \sum _j \lambda _j \ell _{u_j} (Y)^2.  $$
Consequently, using Cauchy-Schwarz inequality,
$$ \sup _{X\in H^2}| H_{ S  , t }F (X)  - F (X) | \leq  \Vert F \Vert
_Q  \left [
 \int _{ \widetilde S}  \widetilde Q(Y) d\mu _{\widetilde S , t }
 (Y)\right ]^{1/2}.   $$
Besides,
$$  \int _{ \widetilde S}  \widetilde Q(Y) d\mu _{\widetilde S , t }
(Y) =
t \sum _j \lambda _j |\pi _S u_j |^2  = t {\rm Tr } (A \pi _S).$$
One proceeds similarly   for the derivatives of $F$ and  one then deduces Proposition \ref{p3.5}.

\hfill$\Box$

We call polynomial function on $B$, any linear combination of functions written as
$x \mapsto (x\cdot a_1 ) \dots  (x \cdot a_m)$ with $a_j\in B'$.

\begin{prop}\label{p3.6} Fix $F\in S(H^2 , Q_A)$. For any finite dimensional $E$ with $E\subset B' \subset H$ there exists a bounded operator $T_E (F)$ in $L^2 ( B, \mu _{B ,
h/2})$  satisfying for   any  polynomials functions $f$ and $g$ on
$B$,
\bl{3.15}  Q_h ^{weyl } ( H_{ ( E^{\perp})^ 2  , h/2 } F) (f , g) =  < T_E
(F) f , g>.\el
 Moreover, the following estimate holds,
\bl{3.16}\Vert T_E (F) \Vert \leq  \Vert F \Vert _Q \
D(E^2 , q_A  |_{E^2 } , 81 \pi h
\max (1, \Vert  \pi _E {\cal F}  A i_E  \Vert _{E^2 , q_A|_{E^2} }  )
, \pi _E {\cal F}  A i_E) ^{1/2}.\el
\end{prop}

{\it Proof.} \newline
Since $f$ and $g$ are polynomials functions, these two functions are are linear combinations of functions $x \mapsto (x\cdot a_1 ) \dots  (x \cdot a_m)$ with $a_j$ in $B'$. There is a finite dimension subspace $S$ being orthogonal to $E$ and such that, for all $j$, $a_j$ belongs to $E \oplus S$. One sets $B_1 = E$, $B_2 = S$ and set  $B_3=\{x\in B\ |\ u(x) = 0,\ \forall u\in E \oplus
S\}$. One has $B = B_1 \oplus B_2 \oplus B_3$. The variables $X$ in $B^2$ may be written as $ X = (X' ,
X'' , X''')$ with $ X'\in B_1^2$, $X''\in B_2^2$,
$X''' \in B_3^2$. Using Definition \ref{d2.6}, Proposition \ref{p3.3} and  $H_{ ( E^{\perp})^ 2 , h/2} =
H_{ S^ 2 , h/2} H_{ B_3^ 2 , h/2} $ and noticing that $ H _h^{gauss} (f
, g )$ is a cylindrical function of basis $(E \oplus S)^2$, one has,
$$ Q_h ^{weyl } ( H_{ ( E^{\perp})^ 2  , h/2 } F) (f , g) =$$
$$\int_{B^2  }
 \widetilde (H_{S^2 , h/2} H_{B_3^2 , h/2}F ( X' , X''  , X''' )
\ \widetilde H _h^{gauss} (f , g ) ( X', X'')   d\mu _{ B^2, h/2} (X', X'' , X'''). $$
Since $E$ and $S$ are finite dimensional subspaces, the stochastic extension is concerned only with the variables $X''$. In view of the heat operator semigroup property, one deduces
$$ \int \widetilde H_{B_3^2 , h/2} F ( X' , X''  , X''' )
d\mu _{B_3^2 , h/2 } (X''') = \int \widetilde F ( X' , X''  , X''')
d\mu _{B_3^2 , h } (X''').$$
This function is denoted by $(Av _{E
\oplus S, h} F ) (X' ,  X'')$ and it belongs to $S(H^2, Q)$. Consequently,
$$ Q_h ^{weyl } ( H_{ ( E^{\perp})^ 2  , h/2 } F) (f , g) =
\int_{ (E\oplus S)^2}  (H_{S , h/2} Av _{E\oplus S , h} F ) (X) \
\widetilde H _h^{gauss} (f , g ) ( X)
 d\mu _{ (E\oplus S)^2, h/2} (X). $$
Applying Proposition \ref{p3.4} with
$H=E\oplus S$ and $\Phi =  Av _{E\oplus S , h} F$, one has
$$  Q_h ^{weyl } ( H_{ ( E^{\perp})^ 2  , h/2 } F) (f , g) =
\int _{S^2} Q_h ^{weyl, E} (Av _{E\oplus S , h} F )( \cdot , X'') )
   ( T_{  X'', h }f ,  T_{  X'', h }g) d\mu _{ S^2, h} (X''). $$
Applying Theorem \ref{t1.4} with the space $E$ and the quadratic form $Q_{A , E}$ being the restriction of $Q_A$ to $E^2$ together with the
function $G(X') = (Av _{E\oplus S , h} F )( X', X'')$ satisfying
$\Vert G  \Vert _{Q_{A , E}}  \leq \Vert F \Vert _{Q_A}$, one obtains
$$ | Q_h ^{weyl } ( H_{ ( E^{\perp})^ 2  , h/2 } F) (f , g) | \leq
\Vert F \Vert _{Q_A} \ D(E^2 , 81 \pi h S(A,E)  , Q_{A , E}) ^{1/2}
\ \dots  $$ 
$$ \dots  \int _{S^2}
\Vert T_{  X'', h }f \Vert _{L^2 (E , \mu   _{ E , h/2} )}\ \Vert T_{
X'' , h} g \Vert _{L^2 (E , \mu   _{ E , h/2})}  \ d\mu _{ S^2, h}
(X''), $$ 
where $S(A,E)=\sup (1,||\pi_E{\cal F}A i_E||_{E^2,q_A|_{E^2}})
$. The Segal Bargmann transform being unitary, one has for any $f$ polynomial function $f$,
$$ \int _{S^2}
\Vert T_{  X'', h } f \Vert _{L^2 (E , \mu   _{ E , h/2} )}^2 \
 d\mu _{ S^2, h} (X'') = \Vert f \Vert ^2.$$
Thus,  there exists an operator $T_E (F)$ being
bounded in $L^2 ( B , \mu _{B , h/2})$  and satisfying (\ref{3.15}) and
 (\ref{3.16}). The proof of Proposition \ref{p3.6} is complete.

 \hfill$\Box$

{\it Proof of Theorem \ref{t1.6}. } \newline
 Let $(E_n)$ be an increasing sequence of finite dimensional subspaces with dense union  in $H$.

{\it Step 1.} One applies Proposition \ref{p3.6} with  $E=E_n$. Thus, one obtains a sequence of operators $T_{E_n}F$ satisfying  (\ref{3.15})
and (\ref{3.16}). One notices that, for any finite dimension subspace $E$ of $H$ and for all$z\in {\bf C}$,
 $$ D(E^2 , q_A |_{E^2 } , z , \pi
 _E {\cal F}  A i_E)  =
 D( B_A , q_A , z  , \pi _E {\cal F}  A P_E )$$
and
 $$  \Vert  \pi _E {\cal F}  A i_E  \Vert _{E^2 , q_A|_{E^2} } =
  \Vert \pi _E {\cal F}  A P_E  \Vert _{B_A , q_A}.  $$
From Proposition \ref{p3.2},
  $$ \lim _{n \rightarrow \infty } \Vert   \pi_{E_n} {\cal F} A
  P_{E_n} - {\cal F} A \Vert
  _{{\cal L}^1 ( B_A, q_A) } = 0,$$
implying,
  $$ \lim _{n \rightarrow \infty }  \Vert   \pi_{E_n} {\cal F} A
  P_{E_n}  \Vert _{
  E_n^2 , q_A|_{E_n^2}}  = \Vert {\cal F} A \Vert _{B_A, q_A }. $$
In view of continuity properties concerning the Fredholm determinant and setting,
  $$ D_n =   D(E_n^2 , q_A  |_{E_n^2 } , 81 \pi h
\max (1, \Vert  \pi _{E_n} {\cal F}  A i_{E_n}  \Vert _{E_n^2 ,
q_A|_{E_n^2} }  )  ,
 \pi _{E_n}  {\cal F}  A i_{E_n}), $$
one has
 $$ \lim _{n \rightarrow \infty }  D_n =
 D(B_A , q_A ,  81 \pi h
\max (1, \Vert   {\cal F}  A   \Vert _{B_A , q_A }  ) \, \  {\cal F}
A  ). $$ 

{\it Step 2.} Let us prove here that
 $T_{E_n} (F) $ is a Cauchy sequence. Set $m,n$ with  $m< n$ and define $S_{m,n}$ as the
orthogonal supplement of  $E_m$ in $E_n$. One then has,
$$ H_{ ( E_m^{\perp}) 2  , h/2 } = H_{ ( E_n^{\perp}) 2  , h/2 }
H_{ S_{mn} , h/2 }. $$
Consequently,
$$ T_{E_m} (F)  = T_{E_n} (H_{ S_{mn} , h/2 } F).$$
From Proposition \ref{p3.6} and Proposition \ref{p3.5}, one sees
$$ \Vert T_{E_n} (F)  - T_{E_m} (F) \Vert \leq \Vert ( I - H_{ S_{mn}
, h/2 } ) F
\Vert _Q  D_n \leq (h/2)^{1/2} D_n \Vert  F \Vert _Q   {\rm Tr} (A
\pi _{S_{mn} } ) ^{1/2}. $$ 
The sequence $D_n$ is bounded according to Step 1 and the trace in the right hand side tends to zero when $m\rightarrow\infty$  in view of Theorem 6.3 in \cite{GK}. This implies that $T_{E_n} (F) $ is a
Cauchy sequence in ${\cal L} ( L^2 ( B , \mu _{B , h/2}))$.

{\it Step 3.}  From step 2, it is known that $\lim_{n\rightarrow\infty}T_{E_n}(F)=T$  exists and belongs to ${\cal L} ( L^2 ( B , \mu _{B , h/2}))$. This limit $T$ is denoted by
$Op_h^{weyl } (F)$. We now derive identity (\ref{1.6}) for any polynomial functions $f$ and $g$. To this end,  $n$ goes to $\infty$ in equality (\ref{3.15}). The right hand side tends to $< Tf , g>= < Op_h ^{weyl}(F) f , g>$. From (\ref{2.7}), for $f$ and $g$ being (stochastic extensions of) polynomials functions,
$$ | Q_h ^{weyl } ( F - H_{ ( E_n^{\perp}) 2  , h/2 } F) (f , g) |
\leq  \Vert  \widetilde F - H_{ ( E_n^{\perp}) 2  , h/2 } \widetilde
F \Vert _{\infty}  C(f , g),$$
where $C(f , g)$ is defined in (\ref{2.7}). In the above right hand side, the norm is the $L^\infty$ norm on $B^2$. From Proposition \ref{p3.5}, one has
$$ \Vert  \widetilde F - H_{ ( E_n^{\perp}) 2  , h/2 }
\widetilde  F \Vert _{\infty}
 \leq
\Vert   F - H_{ ( E_n^{\perp}) 2  , h/2 } F \Vert _{Q} \leq \Vert   F
\Vert _{Q} (h/2)^{1/2} {\rm Tr } (A \pi _{(E_n^{\perp}) 2 } ) ^{1/2}.
$$
The trace term in the above right hand side tends to zero as $n$ goes to $\infty$ using Theorem  6.3 in \cite{GK}. Consequently, the left hand side of
 (\ref{3.15}) tends to
 $Q_h^{weyl} (F) (f , g)$ and the identity (\ref{1.6}) is proved.

 \hfill$\Box$

\section{Reduced metaplectic group covariance.}\label{s4}

\subsection{Introduction.}\label{s4.1}

The symplectic form on $H^2$ is defined by,
\bl{4.1} \sigma ( (x , \xi), (y, \eta)) = y\cdot \xi - x\cdot \eta.
\el
A linear mapping keeping invariant the symplectic form is symplectic. In Shale \cite{Sha} (in the Fock spaces framework) or in B.
Lascar \cite{Las}  (with the $L^2$ setting), one associates unitary operators $U_{\chi}$ in ${\cal F}_s (H_{\bf C})$ or in $L^2(B , \mu _{B , h/2})$ with suitable symplectic maps $\chi$ in $H^2$. This extends  metaplectic transforms to the infinite dimension case. These operators  are also called  Bogoliubov transforms. Let us begin with recalling the definition of the suitable symplectic maps.

\begin{defi}\label{d4.1}  We denote by ${\rm Sp}_2 (H^2)$ (resp.
${\rm Sp}_1 (H^2)$)
the group of symplectic maps  satisfying that  $ [\chi , {\cal
F}]$ is an Hilbert-Schmidt operator, (resp. a trace class operator) where ${\cal F} (x , \xi ) = (-\xi, x)$.
\end{defi}

The fact that  $\chi$ is symplectic can be written as $^t  \chi {\cal
F} \chi = {\cal F}$. Then, a symplectic map $\chi$ belongs to ${\rm Sp}_2
(H^2)$ (resp. ${\rm Sp}_1 (H^2)$ ) if and only if
 $I - ^t \chi \chi $ is an Hilbert Schmidt operator (resp. trace class operator).
The following Theorem is a standard result of  Shale \cite{Sha}.

\begin{theo}\label{t4.2} For any $\chi $ in  $ {\rm
Sp}_2 (H)$ there exists an unitary operator  $U_{\chi}$ in $L^2 (B , \mu _{B , h/2})$  satisfying  for any $X\in H^2$,
\bl{4.2} U_{\chi} ^{\star} V_h (X)  U_{\chi} = V_h ( \chi^{-1}  (X) ),\el
where Weyl translations $V_h(X)$ are defined in (\ref{2.2}).
For  any $\chi_1$ and $\chi_2$ in ${\rm Sp}_2 (H)$, one has 
\bl{4.3} U_{\chi_1} \circ U_{\chi_2}= C( \chi_1 , \chi_2) U _{\chi_1
\chi_2},\el where $ C( \chi_1 , \chi_2)$ is a constant of module 1.
\end{theo}

An explicit form of the  $U_{\chi}$ can be found in \cite{BD} in the Fock space framework and in \cite{Las} in  $L^2$ space setting.

\subsection{Beals property covariance.}\label{s4.2}

When $F$ is a real linear form on $H^2$ one may define (see Definition \ref{d2.6}) a quadratic form
$Q_h^{weyl } (F)$ on the space of polynomial functions stochastic extensions. It is the quadratic form corresponding to the selfadjoint unbounded operator denoted by $Op_h^{weyl} (F)$ with domain containing polynomial functions stochastic extensions. In the Fock space framework, this operator is usually called Segal field.

From Proposition 8.10 of \cite{AJN-JFA}, if $ F(x , \xi) = a\cdot x + b\cdot \xi $ with $(a , b)\in H^2$ then $Op_h ^{weyl } (F)$ is the following selfadjoint operator defined on a suitable domain,
\bl{4.4} Op_h ^{weyl } (F) f(u) = \ell _{ a+ib} (u) f(u) + {h\over i} b \cdot {\partial f
\over \partial u}.\el

Thus, the operators $V_h
(X)$ in Section \ref{s2} and the operator $L_h (X)$ in Definition
\ref{d1.7} are related by,
\bl{4.5} V_h (X) = e^{-{i\over h} L_h (X)},\quad X\in H^2.\el
From Theorem \ref{t4.2}, one deduces that, if $F$ is a continuous linear form on $H^2$,
 \bl{4.6}U_{\chi} ^{\star} Op_h ^{weyl } (F)  U_{\chi} = Op_h ^{weyl } (F\circ \chi).\el
In particular, for any $V\in H^2$
\bl{4.7}
U_{\chi} ^{\star} L_h(V)  U_{\chi} = L_h ( \chi^{-1} (V)).\el

\begin{prop}\label{p4.3} Let $\chi \in {\rm Sp}_2 ( H^2)$ and
 $U_{\chi}$ be the unitary map in Theorem \ref{t4.2}. Then, the  map $A \mapsto
 U_{\chi}^{\star }
 A U_{\chi}$  is an isometry from the space $\in  {\cal L} (Q)$ (see Definition \ref{d1.7}) onto the space
 $ {\cal L}_p (Q\circ \chi)$.
\end{prop}

 {\it Proof.} \newline
 Assume $p=1$. Let $(A_h)\in {\cal L}_1 (Q)$. For all $V\in H^2$, one has in view of (\ref{4.7}),
 $$  {\rm ad} ( L_h V) (U_{\chi}^{\star } A_h U_{\chi} ) =
 U_{\chi}^{\star }
   \Big [ (U_{\chi} ( L_h V)U_{\chi}^{\star } ) , A_h \Big ] U_{\chi}
   $$
 $$ =  U_{\chi}^{\star }  \Big [ L_h ( \chi(V)) , A_h \Big ] U_{\chi}
 $$
and using $(A_h)\in {\cal L}_1 (Q)$ and since  $U_{\chi}$ is unitary, one gets that,
 $$ \Vert {\rm ad} ( L_h V) (U_{\chi}^{\star } A_h U_{\chi} ) \Vert
 \leq \Vert A_h \Vert _{{\cal L}_1 (Q) } Q ( \chi (V)) ^{1/2}.
 $$
Consequently, $U_{\chi}^{\star } A_h U_{\chi}$ indeed belongs to $ {\cal L}_1 (Q\circ {\chi})$.  One proceeds similarly for the space ${\cal L}_p (Q)$.

\hfill$\Box$
\subsection{Weyl calculus covariance.}\label{s4.3}

For each $F\in S(H^2 , Q_A)$ and any $\chi \in {\rm
 Sp}_2 (H^2)$, the function
 $F\circ \chi \in S(H^2 , Q_A\circ \chi)$. In particular, these functions $F$ and $F\circ \chi$ have stochastic extensions. Thus, one may consider the  Weyl operators associated to $F$ and $F\circ \chi$.

 \begin{theo}\label{t4.4} Set $\chi \in {\rm
 Sp}_1 (H^2)$,
  $F\in S(H^2 , Q_A)$ where $A$ is a nonnegative selfadjoint trace class operator in $H^2$. Then,
 \bl{4.8} U_{\chi}^{\star} Op_h ^{weyl} (F)  U_{\chi} = Op_h ^{weyl}
 (F\circ \chi),\el
 where $ U_{\chi}$ is the unitary operator in Theorem \ref{t4.2}.
\end{theo}

 Proposition 2.3 in  \cite{Las} with slight modifications reads as,

\begin{theo}\label{t4.5}  Each $\chi \in {\rm Sp}_2 (H)$,
(resp.  $ {\rm Sp}_1 (H)$) can be uniquely written as $ \chi =
\chi_1 \chi_2 \chi_3$ where,

i)  The map $\chi_1 \in {\cal L} (H^2)$  is unitary and symplectic.

ii)   The map $\chi_2 \in {\cal L} (H^2)$ is written as $(x ,
\xi) \mapsto ( x , \xi + Sx)$ where $ S \in {\cal L} (H)$ is a symmetric and Hilbert-Schmidt (resp. trace class) operator.

iii)   The map $\chi_3 \in {\cal L} (H^2)$ is written as $(x ,
\xi) \mapsto  (T x ,  T^{-1}\xi)$ where $ T \in { GL} (H)$
is symmetric positive invertible operator satisfying  $ T  -
I$ is an Hilbert-Schmidt (resp. trace class) operator.
\end{theo}

{\it Proof.}\newline Set $V = \{ 0 \} \times H \subset H^2$,
$\chi \in {\rm Sp}_2 (H)$.  Then $\chi (V)$ is a closed subspace of $H^2$ on which the symplectic form $\sigma$ vanishes identically. Let  $(u_j)$ be an Hilbert basis of $\chi (V)$ (endowed with the usual scalar product in $H^2$) and let $(e_j)$ be an Hilbert basis of
$H$. There exists $\varphi _1$ orthogonal and symplectic in $H^2$,
satisfying  $\varphi_1 ( 0, e_j) = u_j$. Thus, $\varphi = \varphi_1
^{-1} \chi $ maps $V$ into $V$. One then  may write,
$$\varphi = \begin{pmatrix} A & 0 \cr C & D    \end{pmatrix}.$$
One knows that a symplectic map is unitary if and only if it commutes with ${\cal F}$. Consequently,
$\varphi$ belongs to ${\rm Sp}_2 (H)$ (resp. ${\rm Sp}_1 (H)$).
Thus, $A$ is invertible, $D = ^t A^{-1}$, $A-D$ and $C$
are Hilbert-Schmidt operators. Using polar decomposition one writes  $A =
US$ with $U$  orthogonal and $S$ symmetric, positive  and invertible. One has $D = U S^{-1}$. Therefore, $S - S^{-1}$ is an
Hilbert Schmidt (resp. trace class) operator implying that $S-I$
Hilbert Schmidt (resp. trace class). Consequently,
$$  \varphi =  \psi _1 \chi_2 \chi_3, \quad
\psi_1 =
\begin{pmatrix}  U & 0 \cr 0 &  U   \end{pmatrix},  \quad
 \psi_2 =
\begin{pmatrix}  I & 0 \cr T &  I    \end{pmatrix},  \quad
 \psi_3 =
\begin{pmatrix}  S & 0 \cr 0 &  S^{-1}    \end{pmatrix}, $$
 where $T = U^{-1} C S^{-1} $. One has $
\chi = \chi_1 \chi_2 \chi_3$ with $\chi_1 = \varphi_1 \psi_1$ which is symplectic and unitary. Since $\psi_2$ is symplectic, then
$T$ is an Hilbert Schmidt (resp. trace class) operators.
Therefore, $\chi_2$ and $\chi_3$ have the required properties.

\hfill $\Box$

{\it Metaplectic transforms explicit forms.} \newline
We recall here, without proof the explicit form of  $U_{\chi}$ given by B. Lascar \cite{Las}
in each case i), ii) and iii) above (Theorem \ref{t4.5}).

{\it Case of $\chi_1$.} The Bargmann transform (also named FBI transform) of $f$ is the function $T_h f$ defined on $H^2$ by,
\bl{4.9}
f\mapsto
 (T_hf) (X) = { < f , \Psi_{Xh} > \over  < \Psi_{0h} , \Psi_{Xh}
>},\quad X\in H^2.\el
We know from \cite{AJN-JFA}  that the map $X \mapsto (T_hf) (X)$ defined on $H^2$ has a stochastic extension denoted by  $\widetilde T_hf$ in  $L^2( B^2, \mu _{B^2 , h})$. In addition, the map $f\mapsto \widetilde T_h f$ is a partial isometry from $L^2 (B , \mu _{B , h/2})$
into $L^2 (B^2 , \mu _{B , h})$. For any unitary and symplectic $\chi \in {\cal L} (H^2)$ one associates as in  \cite{Las} an operator $U_\chi$ defined by,
\bl{4.10}  U_{\chi} f = \widetilde T_h ^{\star} ( (\widetilde T_hf) \circ
\chi ).\el The above  right hand side is the stochastic extension
of the map defined on $H^2$ by $X \mapsto (T_hf) ( \chi (X))$.

{\it Case of  $\chi_2$.} Let $(u_j)$ be an Hilbert basis of $H$
constituted with eigenvectors of $S$ and let  $\lambda _j $ be the corresponding eigenvalues. We set,
\bl{4.11} \varphi_S(x) = \sum _j \lambda _j ( \ell _{u_j} (x) ^2 -{h\over 2} ).\el
The above series defines an element of  $L^2(B , \mu _{B , h/2})$.
Indeed, it is known that for all $u\in H$, the element
$u\otimes u$ in ${\cal F}_s (H_{\bf C})$ corresponds through the Segal isomorphism to the element  $ \ell _u(x)^2 - (h/2)$ in $L^2(B , \mu_{B , h/2})$ (\cite{Jan}\cite{Sim-PP2}). The series $\sum \lambda _j u_j \otimes u_j$ thus defines an element in ${\cal F}_s (H_{\bf C})$ corresponding to $\varphi_S$ through this isomorphism. The metaplectic transform related to $\chi_2$, $U_{\chi_2}$, is defined in \cite{Las} by,
\bl{4.12} U_{\chi_2} f (x)  = e^{ {i\over 2h} \varphi_S(x) } f(x). \el
For  $\chi_3$,
we don't follow \cite{Las}  but we start with two results concerning change of variables in integrals on Wiener spaces. They will be used again in the proof of Theorem \ref{t4.4}.

\begin{prop}\label{p4.6} Let  $(i, H, B)$ be
 a Wiener space. Set $\mu _{B , h}$ the Gaussian measure on $B$ with variance $h$. Set $\varphi \in GL (B)$ and let  $\nu_h$ be the Borelian measure on $B$ defined by,
$$ \nu_h (E)= \mu _{B , h} (  \varphi ^{-1} (E) )\ {\rm with}\ E\subset B,\ E {\rm \ a\ Borel\ set}.$$
If  $\varphi $ restricted to $H$ satisfies $\varphi ^t\circ \varphi = I + K$ where $K\in {\cal L}(H)$
is selfadjoint and Hilbert Schmidt then $\nu_h$ is absolutely continuous with respect to $\mu _{B, h}$. Set
$D_h (x)$ the density of $\nu_h$ with respect to $\mu _{B , h}$.
Then, this function $D_h$ belongs to $L^1 (B , \mu_{B , h} )$ and one has, for every $\Phi \in L^1 ( B , \mu _{B , h})$,
\bl{4.13} \int _{B}  \Phi  ( x) D_h (x) d\mu _{B , h} (x) =
 \int _{B}  \Phi ( \varphi (x) )  d\mu _{B , h} (x). \el
In addition, if $(u_j)$ is an Hilbert basis of $H$ made of eigenvectors of $K$ and if the $\nu _j$ are the related eigenvalues, then one has,
\bl{4.14} D_h (x) = \prod (\nu_j +1)^{-1/2} \exp \Big [ \nu _j \ell_{u_j}
 (x)^2 /
 ( h ( \nu_j +1))\Big  ].\el
\end{prop}

{\it Proof.}\newline
 For any  measurable and bounded $\Phi$ one has,
$$ \int _{B} \Phi ( x) d\nu_h (x) =
 \int _{B}  \Phi ( \varphi   (x) )  d\mu _{B , h} (x).$$
Thus, for any  $a\in B'$, 
$$ \int _{B}  e^{-i <a , x>} d\nu_h (x) =
 \int _{B}   e^{-i <a , \varphi (  x ) >}   d\mu _{B , h} (x).$$
From Theorem \ref{t1.2},
$$ \int _{B}  e^{-i <a , x>} d\nu_h (x) =
 \int _{B}   e^{-i <^t\varphi a ,    x>}   d\mu _{B, h} (x)
 = e^{-{h\over 2} |^t\varphi a|^2 }.  $$
The dual quadratic form $Q(a) = |^t\varphi a|^2 $ is written as $Q(a)
= (a + Ka) \cdot a$ where $K = \varphi ^t \varphi -I$. We have assumed that $K$ is a Hilbert Schmidt operator. From a result of Segal and Feldman (a statement can be found in Ramer \cite{Ram}, Section 3 or in Theorem 1.23 of \cite {Sim-PP2}), we know that $\nu_h$ is absolutely continuous with respect to $\mu _{B,
h}$. If $D_h (x)$ denotes the corresponding density, one indeed has  (\ref{4.13}). Equality (\ref{4.14}) also is a statement of  Ramer \cite{Ram} (Section 3).

\hfill$\Box$

\begin{prop}\label{p4.7} i) Let $ T \in { GL} (H)$
be a  symmetric positive invertible operator such that  $ T  -
I$ is an Hilbert-Schmidt operator. Then, there exists an
 Hilbert space $B_{T}$ containing $H$   such that the triplet
 $(i, H , B_{T})$ is a Wiener space, and the map $T \in {\cal L} (H)$ has,
  by density, an extension as an element of $GL(B_{T})$.

ii)  Let $\psi$ be one of the symplectic map $\chi_2$ or $\chi_3$ in Theorem \ref{t4.5}. There exists an
 Hilbert space $B_{\psi}$ containing $H^2$   such that the triplet $(i, H^2 , B_{\psi})$ is a Wiener space,
 and such that  the map $\psi \in {\cal L} (H^2)$ has an extension as an element of $GL(B_{\psi})$.
\end{prop}

{\it Proof.}\newline
i) Let  $H_0={\rm Ker}\,T-I$, $H_1=H_0^\perp$
its orthogonal in $H$ and  $B_1$  the completion of $H_1$ for the norm $| (T-I) x|$.
By \cite{Kuo} (Exercise 17, page 59), we know that $(i_1, H_1, B_1)$ is a Wiener space, and  therefore,
if $B_T = H_0 \otimes B_1$, with the norm sum of those of the two components,
then $(i, H, B_T)$ is a Wiener space. If $x = x'+ x''\in H$,
with $x'\in H_0$ and $x''\in H_1$, we have $ |Tx |_{B_T} \leq \max (1, \Vert T \Vert )
|x |_{B_T} $, and therefore $T$ has an extension in ${\cal L} (B_T)$. Since $T^{-1}$
satisfies the same hypotheses, ( $T^{-1} - I$ is Hilbert Schmidt), and since
${\rm Ker}\,T^{-1}-I = {\rm Ker}\,T-I$, $T^{-1}$ has also an extension in
${\cal L} (B_T)$, and the extension of $T^{-1}$ is the inverse of the extension of $T$.

ii) Let $\psi (x, \xi) =( x , \xi + Sx)$ where $ S \in {\cal L} (H)$ is a symmetric and Hilbert-Schmidt
operator. One has $H^2 = {\rm ker } S \oplus ( {\rm ker } S )^{\perp} \oplus H$. Let $B_1$ and
$B_3$  be Hilbert spaces  such that  $(i_1 , {\rm ker } S , B_1)$ and $(i_3 , H , B_3)$ are Wiener spaces. Moreover let $B_2$ be the completion of $( {\rm ker } S )^{\perp}$ with respect to the norm $|Sx |$. Since  $S$ is Hilbert Schmidt  $(i_2 , ( {\rm ker } S )^{\perp} , B_2)$ is a
 Wiener space. One sets $B_{\psi} = B_1 \oplus B_2 \oplus B_3$ with its natural norm. Then $(i, H^2, B_{\psi})$ is a Wiener triplet. For $(x , \xi) = ( x_1 , x_2 , \xi ) $ in $H^2$ one has
$$ |\chi_2 (x , \xi) |_{B_{\psi}} = |x_1| _{B_1} + |Sx_2| + |\xi + S
x_2 |_{B_3}
\leq  |x_1| _{B_1} + |Sx_2| + |\xi|_{B_3}   + C |S x_2 |,$$  where
$C>0$ is  such that  $|\xi |_{B_3} \leq C |\xi |$. Clearly, $ |\chi_2 (x , \xi) |_{B_{\chi_2}}  \leq (1+C) |(x ,
\xi) |_{B_{\chi_2}}$ and there exists an extension of $\chi_2$ in ${\cal L}(B_{\chi_2})$.
We conclude as above that this extension is invertible. For $\psi = \chi_3$ one may proceed in the same way.

\hfill$\Box$

{\it Case of $\chi_3$.} Suppose that $B_1$ and $B_2$ are Banach spaces  such that  $(i_1 , H , B_1)$ and $(i_2 , H, B_2)$ are Wiener spaces. The $L^2(B_j , \mu _{B_j , h/2})$ are isomorphic. Therefore, we are allowed to choose  a convenient Banach space $B$ adapted to $\chi_3$. If $\chi_3 (x , \xi) =  (T x ,  T^{-1}\xi)$, where $ T \in { GL} (H)$
is symmetric positive invertible operator such that  $ T  -I$ is Hilbert-Schmidt. We use the space $B_T$
in the point i) of  Proposition \ref{p4.7}.
Applying Proposition \ref{p4.6} with $B_T$, with the extension $\widetilde T$ of $T$,
with the variance   $h/2$,  and setting,
\bl{4.15} (U_{\chi_3} f ) (x)  = f( \widetilde T x) D_h(\widetilde T
x)^{-1/2},\el
one defines a unitary map in $L^2 (B_T , \mu _{B_T , h/2})$ which can be carried to $L^2 (B ,\mu _{B , h/2})$.

{\it Proof of  Theorem \ref{t4.4}.}
\newline
It is sufficient to prove  (\ref{4.8}) for each kind of generator $\chi_1$,  $\chi_2$ and $\chi_3$ of Theorem \ref{t4.5}. The form of $U_{\chi}$ obtained above for the three generators shows that $U_{\chi}$ commutes with the parity operator $f\mapsto \check f$. Therefore, for any  $f,g\in L^2 (B , \mu _{B , h/2})$ and $X\in H^2$ one has
\bl{4.16} H_h^{gauss} (U _{\chi} f , U _{\chi} g) (X) = e^{ (|X|^2 - |\chi
^{-1} (X)|^2 )/h}
 H_h^{gauss} ( f , g) ( \chi ^{-1} (X)).\el
Therefore the proof of (\ref{4.8}) is only a problem of change of variables in an integral.

{\it Proof of (\ref{4.8}) for $\chi_1$.}\newline
Let $f$ and $g$ be two polynomials functions. One has,
$$ <  U_{\chi_1}^{\star} Op_h ^{weyl} (F)  U_{\chi_1} f , g> = \int
_{B^2}
\widetilde F (X) \widetilde H_h^{gauss} (U _{\chi_1} f , U _{\chi_1}
g) (X) d\mu _{B^2 , h/2} (X). $$
It follows from (\ref{4.16})  that $H_h^{gauss} (U _{\chi_1} f , U _{\chi_1} g) (X)$ has a stochastic extension and since $\chi_1$ is unitary,
$$H_h^{gauss} (U _{\chi_1} f , U _{\chi_1} g) (X) =
 H_h^{gauss} ( f , g) ( \chi_1 ^{-1} (X)). $$
Consequently, the left hand side is a polynomial function which has therefore a stochastic extension.
Since $F$ is in $S(H^2 , Q_A)$, it also has a stochastic extension.
 Let
$\Phi (X) = F(X) H_h^{gauss} ( f , g) ( \chi_1 ^{-1} (X))$, and $\widetilde \Phi$ its stochastic extension in $L^1 (B^2 , \mu _{B^2 , h/2})$. Let
$\widetilde \Phi_{\chi}$ be the extension  of $\Phi\circ {\chi}$ (which exists for the same reason).
We shall prove that if  $\chi$ is unitary,
$$
\int_{B^2} \widetilde \Phi_{\chi}(X)d\mu_{B^2,h/2}(X) = \int_{B^2}
\widetilde \Phi(X)d\mu_{B^2,h/2}(X).
$$
In fact, if $(E_n)$ is an increasing sequence of finite dimensional subspaces with dense union in  $H^2$ one has,
$$
\int_{B^2} (\Phi \circ\chi)(\widetilde{\pi}_{E_n}(X))d\mu_{B^2,h}(X)
= \int_{E_n} (\Phi \circ\chi)(x)d\mu_{E_n,h}(X) $$
$$
= \int_{\chi E_n} \Phi (x)d\mu_{\chi E_n,h}(X) = \int_{B^2}
\Phi(\tilde{\pi}_{\chi E_n}(X))d\mu_{B^2,h}(X). $$
Letting $n\rightarrow \infty$ one concludes easily using Lebesgue theorem. Thus, if $f$ and $g$ are polynomials functions,
$$ <  U_{\chi_1}^{\star} Op_h ^{weyl} (F)  U_{\chi_1} f , g> =
<   Op_h ^{weyl} (F\circ \chi_1 ) f , g >. $$
For general  $f$ and $g$, the  equality is obtained by density. Consequently,
 (\ref{4.8})  is proved for $\chi=\chi_1$.

Before the proof of (\ref{4.8}) for $\chi_2$ and $\chi_3$, we need the additional two  propositions.

\begin{prop}\label{p4.8}   Let $(i, H, B)$ be a Wiener space,  $T$ a trace class selfadjoint
operator  in $H$ and let the $\lambda _j $ be the eigenvalues of $T$. If
 $p \lambda _j  <1$ for all $j$, then the function  $F$ defined on $H$ by  $F(X) =
e^{ (TX)\cdot X/h}$ has a stochastic extension  to $L^p ( B ,
\mu _ { B , h/2})$.
\end{prop}

{\it Proof.}\newline Let $(u_j)$ be an Hilbert basis of $H$ consisting of eigenvectors of $T$, let  $p'$ with $p < p' $ and $p' \sup
{\lambda _j} < 1$ and let $q' =  p p' / ( p' - p)$. It is known that $\varphi (X) = (TX) \cdot X$
has a stochastic extension stochastic in $L^{q'} ( B , \mu _{B , h/2})$, which is given by,
$$ \widetilde \varphi (X) = \sum \lambda _j \ell _{u_j} (x)^2.$$
 It follows that, for each increasing sequence  $(E_n)$ of finite dimensional subspaces of $H$ with dense union,  the sequence $\varphi \circ \widetilde \pi _{E_n}$ is a
Cauchy sequence in  $L^{q'} ( B , \mu _{B , h/2})$.  Let us also prove that the sequence  $ e^{ (1/h) \varphi \circ \widetilde \pi _{E_n}}$ is bounded in  $L^{p'} ( B , \mu _{B , h/2})$.  Indeed,
$$ \int _B e^{ (p'/h) \varphi \circ \widetilde \pi _{E_n}(X)}
d\mu _{B , h/2} (X) =   {\rm det }( I -  p' \pi_{E_n}  T\pi_{E_n}  )
^{-1/2}.$$
  Since $p' \sup {\lambda _j} < 1$,  the integral is well defined and bounded when $ n\rightarrow \infty$.
Since the sequence   $\varphi \circ \widetilde \pi _{E_n}$ is a
Cauchy  sequence in  $L^{q'} ( B , \mu _{B , h/2})$ and since the sequence $ e^{ (1/h) \varphi \circ \widetilde \pi _{E_n}}$ is bounded in  $L^{p'} ( B , \mu _{B , h/2})$, then the  H\"older inequality proves that the sequence
$e^{ (1/h) \varphi \circ \widetilde \pi _{E_n}}$  is a  Cauchy sequence in $L^{p} ( B , \mu _{B , h/2})$,
which proves the proposition.

\hfill$\Box$

\begin{prop}\label{p4.9} Let $f$ and $g$ be
polynomial functions,  let $\psi$ be one of the maps $\chi_2$ or $\chi_3$
  in Theorem \ref{t4.5}. Let $B_{\psi}$ be the Hilbert space obtained in Proposition \ref{p4.7} and we also call $\psi$ the density extension of $\psi$. Set $U _{\psi }$ the metaplectic transform in Theorem \ref{t4.2}. Then, the Wigner function $H_h^{gauss} (U _{\psi } f, U _{\psi } g)$ defined on $H^2$ has a stochastic extension in some $L^p(B_{\psi},\mu_{B_{\psi},h/2})$. This stochastic extension equals to,
\bl{4.17} \widetilde H_h^{gauss} (U _{\psi } f, U _{\psi } g) (X) =
D_h (X) \widetilde H_h^{gauss} ( f,  g) (\psi ^{-1} (X)),\el
 where $D_h (X)$ is the Radon Nikodym density of
Proposition \ref{p4.6}.
\end{prop}

 {\it Proof.}\newline
One uses equality (\ref{4.16}). Since the function $F_1 (X) = H_h^{gauss} ( f , g) ( \chi_j ^{-1} (X))$ is polynomial then
it has a stochastic extension, for $j= 2$ and
 $j= 3$. We may write $$F_2 (X) = e^{ (|X|^2 - |\chi_j ^{-1} (X)|^2 )/h} = e^{ (TX)\cdot
 X/h},$$
 with $T = I - ^t \chi_j ^{-1} \chi_j ^{-1} $. When $j =2$,
(see Theorem \ref{t4.5}),
 $$ \chi_2 = \begin{pmatrix}  I & 0 \cr  S & I  \end{pmatrix},\quad T =
  \begin{pmatrix} -S^2 &  S \cr  S & 0   \end{pmatrix},$$
 under the assumptions of  Theorem \ref{t4.5}. If $\chi_2\in Sp_1(H^2)$ then $T$ is trace class. Let $(v_j)$ be an Hilbert basis of $H$ consisting of eigenvectors of  $S$ and let $(\mu_j)$ be the corresponding eigenvalues. The eigenvalues of $T$ are, $$\lambda _j^{\pm} = {1\over 2}
\Big [ - \mu_j^2 \pm \sqrt { \mu_j ^4 + 4 \mu_j^2} \Big ].$$ Let
 $\rho = \sup |\mu_j|.$  If $p$ satisfies,
$ 1 \leq p < {1\over 2} \left [ 1 + \sqrt { 1 + {4 \over \rho ^2 }}
\right]
$ then, for any $j$, one has $  p \lambda _j ^{\pm} <1$. According to Proposition \ref{p4.8}, $F_2$ has a stochastic extension
in $L^p ( B_{\chi_2} , \mu _ { B_{\chi_2} , h/2})$.
Let  $( U_j ^{\pm})$ be an Hilbert basis of $H^2$ consisting of  eigenvectors of $T$, corresponding to the eigenvalues $\lambda _j ^{\pm}$. Then the stochastic extension of $F_2$ is, $$\widetilde F_2(X) = \exp \left [ \sum _j \sum _{\sigma = \pm}
\lambda _j^{\pm} \ell _{ U_j ^{\sigma } } (X)^2 \right ]. $$
One applies Proposition \ref{p4.6} with $B_1 = B_{\chi_2}$
defined in Proposition \ref{p4.7} and with  $\varphi$ the density extension of $\chi_2$ to   $B_1$ according to Proposition \ref{p4.7}. One sees that  $K= \chi _2 ^t \chi _2 -
I$ (in Proposition \ref{p4.6}) and $T= I - ^t \chi_2 ^{-1}
\chi_2 ^{-1} $ are related bt $T = K ( K + I)^{-1}$. Each of the eigenvalues $\lambda _j^{\pm}$ of  $T$ is related to one of the eigenvalues $\nu_k$ of $K$ by
  $\lambda _j^{\pm} = \nu_k / ( 1+ \nu_k)$, with the same eigenvectors. Besides, one has $(1 - \lambda _j ^+)
(1 -\lambda _j ^- ) = 1$. For each $k$ there exists $k'$  such that  $(1+ \nu_k ) (1+ \nu_{k'})=1$. Consequently, the stochastic extension $\widetilde F_2$ is the function $D_h$ of Proposition \ref{p4.7}, which achieves the proof for $\chi_2$. Similar considerations are valid for $\chi_3$.

\hfill$\Box$

{\it Proof of (\ref{4.8}) for $\chi_2$ and $\chi_3$.}\newline
Let $f$ and $g$ be polynomial functions, and $\chi$ be either $\chi_2$, or $\chi_3$. The quadratic form
 $Q_h^{weyl } (U_{\chi} f , U_{\chi} g) $ may be defined with the integral form of Definition \ref{d2.6} under the only condition that the Wigner function  $H_h^{gauss} (U
_{\chi} f , U _{\chi} g)$ has a stochastic extension. Besides, for this stochastic extension, one may replace $B^2$ by any other Banach space $\widetilde{B}$ such that $(i,H^2,\widetilde{B})$ is a Wiener space. One then can use Proposition \ref{p4.7} and the space $B_{\chi }$ obtained there in order to construct this extension. From Proposition \ref{p4.9},  $H_h^{gauss} (U _{\chi_j} f , U _{\chi_j} g)$ has a stochastic
extension in some $L^p ( B_{\chi } , \mu
_{B_{\chi } , h/2})$ and from (\ref{4.17}) one has for  $j= 2,3$, for $f$ and $g$ polynomial functions,
$$  <  U_{\chi_j}^{\star} Op_h ^{weyl} (F)  U_{\chi_j} f , g> = \int
_{B_{\chi _j}}
\widetilde F (X)  D_h (X) \widetilde H_h^{gauss} ( f,  g) (\chi _j
^{-1} (X)) d\mu _{B_{\chi _j} , h/2} (X).  $$
We apply Proposition \ref{p4.6} with, $$ \Phi (X) = \widetilde F (X)\widetilde H_h^{gauss} ( f,  g) (\chi
_j ^{-1} (X)). $$ Proposition \ref{p4.6} thus shows that,
$$  <  U_{\chi_j}^{\star} Op_h ^{weyl} (F)  U_{\chi_j} f , g> = \int
_{B_{\chi _j}}
\widetilde (F \circ \chi) (X) \widetilde H_h^{gauss} ( f,  g)  (X)
d\mu _{B_{\chi _j} , h/2} (X),  $$
then, for polynomial functions, one deduces,
$$<  U_{\chi_j}^{\star} Op_h ^{weyl} (F)  U_{\chi_j} f , g>
 =  <   Op_h ^{weyl} (F\circ \chi )   f , g>.$$
Finally (\ref{4.8}) then follows by density for  $f,g\in L^2 ( B, \mu _{B , h/2})$.

\hfill$\Box$

\section{Beals characterization.}\label{s5}

\subsection{Wick bisymbol.}\label{s5.1}

With families of operators one usually associates various symbols. Even in finite dimension it does not seem possible to deduce the Weyl symbol from the Wick one by an integral relation. However this may be achieved from Wick bisymbol.

Let $B$ be a bounded operator on $L^2(B \mu_{B , h/2})$, let $X\in H^2$ and let $\Psi_{Xh}$ be the coherent state defined in  (\ref{2.8}). The Wick bisymbol is defined by,
\bl{5.1} S_h ( B) (X , Y) = { < B \Psi_{Xh} , \Psi_{Yh} > \over
< \Psi_{X h} , \Psi _{Yh}> } \quad X,Y\in H^2 .\el
The function $S_h ( B) (X , Y) $ is holomorphic in $x+i \xi$, antiholomorphic in $y+i\eta$ and its restriction to the diagonal is the Wick symbol $\sigma_h ^{wick} (B)$.

\begin{prop}\label{p5.1} Let $d_{diag} (S_hB) (X , Y) $  be the derivative of $H^2 \ni V \rightarrow (S_hB) (X+V ,
Y + V) $.  Suppose that $(B_h)$ is a family in ${\cal
L}_m ( Q_A)$ (see Definition \ref{d1.7}) and $V_1$, \dots  ,$V_m\in H^2$. Then, 
\bl{5.2} | d^m_{diag} (S_hB_h) (X , Y) (V_1 , \dots  , V_m) | \leq  \Vert
B_h\Vert
_{{\cal L}_m ( Q_A) } \  e^{{1\over 4h}(|X-Y|^2 )} \ \prod _{j=1}^m
Q_A(V_j) ^{1/2}.\el
\end{prop}

{\it Proof.}\newline
For $X\in H^2$, let  $V_h(X)$  be the Weyl translation (see Section \ref{s2.1}). From (\ref{2.10}),
$$ < B \Psi_{X+ Z, h} , \Psi _{Y+ Z, h}> =  e^{{i\over 2h}
\sigma (X-Y , Z) }
< V_h(Z)^{\star} B V_h(Z) \Psi_{X, h} , \Psi _{Y, h}>. $$ In particular,
$$ <  \Psi_{X+ Z, h} , \Psi _{Y+ Z, h}> = e^{{i\over 2h} \sigma (X-Y
, Z) }
<  \Psi_{X, h} , \Psi _{Y, h}>. $$ Consequently,
$$ (S_hB) (X + Z , Y+Z) = (S_h (V_h(Z)^{\star} B V_h(Z)) ) (X , Y).
$$
Let $L_h Z$ be the operator defined above Definition \ref{d1.7}. From (\ref{4.5}), one has,
 $$ {d \over dt}  V_h(tZ)^{\star} B_h V_h(tZ) \Big|_{t=0}  =
(i/h)  [ L_hZ, B_h ].$$
Consequently,
$$ ( d_{diag} S_hB) (X , Y) (Z) = (i/h)  S_h ( [ L_h Z , B ] ) (X ,
Y).$$
From (\ref{2.11}) and (\ref{5.1}), one has for any $B$ bounded,
$$ | (S_hB) (X , Y) | \leq \Vert B\Vert e^{{1\over 4h}|X-Y|^2 }.$$
Using Definition \ref{d1.7} and iterating, one deduces (\ref{5.2}).

\hfill$\Box$

\subsection{The finite dimensional situation. }\label{s5.2}

\begin{prop}\label{p5.2} Let   $H$ be a finite dimensional Hilbert space, ${\rm dim}\,H=d<\infty$. Let $(e_j )_{1\leq j \leq d }$ be an Hilbert basis of  $H$ and set $m\geq 0$. Then, for any family
$(A_h)$ in ${\cal L}_{m+4d} ( Q_A)$ there exists a function $F_{
h}$ in $S(H^2 , Q_A)$  satisfying,
\bl{5.3} Op_h^{weyl} (F_h) = A_h, \el
\bl{5.4} H_{ h/2} F_{ h}  = \sigma_h^{wick } (A_h),\el
 and for $h<1$,
\bl{5.5} \Vert F_{ h} \Vert _{m, Q_A} \leq  \Vert A_h \Vert _{{\cal L}_{m+
4d }  ( Q_A) }   \prod _{j\leq d}  ( 1  +   h K S^2 \lambda _j) ,\el
where $K$ is some universal constant and,
\bl{5.6} \lambda _j = \max ( Q( e_j , 0), Q(0, e_j) ),\quad
S = \max (1, \sup \lambda _j).\el
\end{prop}

The standard proof of Beals Theorem in finite dimension corresponds to the step 1 below with $I=\{1,\dots,d\}$. However, in the analogue of (\ref{5.5}), one would obtain a constant going surely to infinity with $d$, whereas the below two steps proof allows to get (\ref{5.5}) where, for some cases, the above product can remain bounded as $d$ goes to infinity.

{\it Proof.}\newline
{\it  Step 1.} For  $I\subset\{ 1
, \dots, d\}$, let $E(I)$ be the subspace of $H$ spanned by the
$e_j$  with $j\in I$. Then, for any  $I\subset\{ 1 , \dots,  d\}$, one may construct a function $G_{I, h} $  satisfying,
\bl{5.7} H_{E(I)  , h/2} G_{I , h}  = \sigma_h^{wick } (A_h).\el
The function $G_{I , h}$ is first defined as a distribution by the following relations,
\bl{5.8} G_{I , h} (X) =  \int _{E(I)^4 }
E_h ( X , S_I, T_I)
 K_h ( S_I, T_I ) d\lambda (S_I, T_I),\el
where $\lambda$ is the Lebesgue measure and where,
\bl{5.9} E_h ( X , S_I, T_I) = (S _h (A_h))  \left  (X + S_I + {T_I \over 2},X + S_I  -{T_I \over 2 } \right ) \el
and
\bl{5.10} K_h ( S_I, T_I ) = 2^{|I| } (2 \pi h )^{-2 |I|} e^{ - { 1\over h }|S_I|^2 - { 1\over 4h } |T_I|^2
- { i\over h } \sigma (S_I, T_I) }. \el
In order to prove (\ref{5.7}) at the distribution level, one may see,
\bl{5.11}  H_{E(I)  , h/2} G_{I , h} (X) = \int _{E(I)^4 } H_{h/2}^{S_I} E_h
( X , S_I, T_I)
 \ K_h ( S_I, T_I ) d\lambda (S_I, T_I) \el

 $$=  \int _{E(I)^4 } E_h ( X , S_I, T_I) \
  H_{h/2}^{S_I} K_h ( S_I, T_I ) d\lambda (S_I, T_I), $$
 where $H_{h/2}^{S_I}$ is the heat operator in the
 variable $S_I$.  An explicit computation shows that,
 $$  H_{h/2}^{S_I} K_h ( S_I, T_I ) = (2 \pi h )^{-2 |I|} \exp \left
 [
 {1\over 2h} \left ( S + {T\over 2 } \right ) \cdot
 \left ( \overline S - {\overline T\over 2 } \right ) -
 {1\over 2h}  \left ( \left | S + {T\over 2 } \right |^2 +
 \left | S - {T\over 2 } \right |^2 \right ) \right ].  $$
 Here and contrarily to the notations used elsewhere in this work, $S \cdot \overline T = (s+ i \sigma) \cdot (t-i \tau)$  if $S = (s ,\sigma)$ and $T = (t , \tau)$.
Since  $S _h (A_h)$ is holomorphic with respect to the first variable (when identifying $H^2$ with the complexified of $H$), antiholomorphic with respect to the second one, then the mean value formula proves that the integral (\ref{5.15}) is equal to  $E_h (X , 0, 0)$, i.e, equals to  $\sigma_h^{wick } (A_h)(X)$. Estimates of Proposition \ref{p5.1} and formula (\ref{5.10}) show that one has a convergence difficulty with the variable $T_I$ and that the integral (\ref{5.8}) is meaningful only as a distribution. Integrating by part as in\cite{ALN-BEALS}, one shows that,
\bl{5.12} G_{I , h} (X) =  \int _{E(I)^4 }
 K_h ( S_I, T_I )\  L_I E_h ( X , S_I, T_I)d\lambda (S_I, T_I), \el
where setting $S = (s , \sigma),T = (t , \tau)$,
\bl{5.13} L_I = \prod _{j\in I } \left ( 1 + { t_j^2 \over h} \right )
 ^{-1}
 \left ( 1 + { \tau_j^2 \over h} \right ) ^{-1} \sum _{{\cal
 M}_2(I)}
 a_{\alpha \beta }  ( S_I / \sqrt {h}) h^{(|\alpha | + |\beta |)/2}
 \partial _s ^{\alpha } \partial _{\sigma } ^{\beta }, \el
where ${\cal M}_2(I)$ is the subset of multi-indices $(\alpha ,
 \beta)$  satisfying  $\alpha _j = \beta _j = 0$ if $j\notin I$ and $\alpha _j \leq 2$, $\beta _j \leq 2$
for any $j\in I$. One may find the values of the coefficients $a_{\alpha \beta }  $ in \cite{ALN-BEALS}.
After these transformations,  (\ref{5.8}) is a convergent integral, and this proves that the distribution
  $ G_{I , h}$ defined in (\ref{5.8}) is a continuous and bounded function. From Proposition \ref{p5.1} one has,
\bl{5.14} |\partial _s ^{\alpha } \partial _{\sigma } ^{\beta }
 E_h ( X , S_I, T_I) | \leq   \Vert A_h\Vert
_{{\cal L}_{|\alpha | + |\beta |}  ( Q_A) }  \  e^{{1\over 4h}|T_I|^2
} \ \prod _{j\in I} Q_A(e_j , 0)^{\alpha _j /2}  Q_A(0 , e_j)^{\beta
_j /2}.\el
 Following \cite{ALN-BEALS}, one has,
\bl{5.15}  \int _{E(I)^4 } |  K_h ( S_I, T_I )|  e^{{1\over 4h}|T_I|^2  }
 \
 \prod _{j\in I } \left ( 1 + { t_j^2 \over h} \right ) ^{-1}
 \left ( 1 + { \tau_j^2 \over h} \right ) ^{-1}
 | a_{\alpha \beta }  ( S_I / \sqrt {h}) | d\lambda (S_I , T_I)
 \leq K^{|I|}, \el
where $K$ is a universal constant.

Taking $I = \{ 1, ..., d \}$, we have a proof of the classical Beals theorem in finite
  dimension, but the  estimations of $ |G_{I , h}|$ obtained in this way are in fact too rough
  when $d$ goes to infinity; a much more precise result is obtained in step 2.

{\it Step 2.}  Let $D_j$ be the subspace spanned by the $(e_j ,
0)$ and $(0, e_j)$ and let $H_{D_j , h/2}$ be the corresponding heat operator. For any finite $I\subset\N$, set
 $$ T_{I , h} = \prod _{j\in I } ( I - H_{D_j , h/2}), $$
 and for any $d\in\N$,
\bl{5.16} F_{ h}  =  \sum _{I \subset \{ 1, \dots  , d\}} T_{I , h}
 G_{I , h},\el
 where $G _{I, h}$ is initially defined in (\ref{5.8}) and more precisely given in (\ref{5.12}).

Let us check (\ref{5.4}). One has if $I  \subset \{ 1, \dots  , d\}$,
 $ H_{ h/2} = H_{E(I)  , h/2} H_{E(\{ 1, \dots  , d\} \setminus I)  ,
 h/2} $
 and by (\ref{5.7}) and (\ref{5.16}),
 $$ H_{ h/2} F_{ h} =  \sum _{I \subset \{ 1, \dots  , d\}}
 H_{E(\{ 1, \dots  , d\} \setminus I)  , h/2}  T_{I , h}
   \sigma_h^{wick } (A_h). $$
All these operators are commuting. One directly checks,
 $$  \sum _{I \subset \{ 1, \dots  , d\}}
 H_{E(\{ 1, \dots  , d\} \setminus I)  , h/2}  T_{I , h}  = I.$$
Thus, one obtains (\ref{5.4}).  Let us now check  estimate (\ref{5.5}).

 One has for $I\subset\{ 1 , \dots  d \}$,
\bl{5.17} T_{I , h} G_{I , h} (X) =  \int _{E(I)^4 }
 K_h ( S_I, T_I )\  L_I T_{I , h} E_h ( X , S_I, T_I)d\lambda (S_I,
 T_I), \el
where $T_{I, h}$ acts on the variable $X$ and where $E_h$ and $K_h$
defined in  (\ref{5.9}) and (\ref{5.10}). One proceeds like in \cite{AJN-JFA},
 Lemma 5.4. We can write  $I - H_{D_j , h/2} $ in three different ways,
 $$  I - H_{D_j , h/2}  = A_j =h^{1/2} ( B_j \partial _{x_j} + C_j
 \partial _{\xi _j})
 = h E_j \Delta _j,$$
 where $A_j$, $B_j$, $C_j$ and $E_j$ are bounded operators in $C_b(H^2)$ with a norm bounded by some universal constant $K_2$. In view of  (\ref{5.9}), one may differentiate $E_h( X , S_I, T_I)$ with respect to $X_I $ or to $S_I$.
 One writes for any all $(\alpha , \beta)\in{\cal M}_2(I)$,
 $$ T_{I , h} \partial _x ^{\alpha } \partial _{\xi } ^{\beta }  =
 \prod _{j \in I} U_j \partial _{x_j }^{\alpha _j }
 \partial _{\xi_j }^{\beta _j },
  $$
where
 $$ U_j  = \left \{ \begin{matrix} A_j \hfill  & {\rm if} & \alpha _j +
 \beta _j \geq 2 \cr
 h^{1/2} ( B_j \partial _{x_j} + C_j \partial _{\xi _j}) & {\rm if}
 &
 \alpha _j + \beta _j = 1 \cr
 h E_j \Delta _j  \hfill  & {\rm if} & \alpha _j + \beta _j  = 0 \cr
 \end{matrix}
 \right .. $$
Consequently,  for any multi-indices $(\alpha , \beta)\in {\cal M}_2(I)$,
\bl{5.18}h^{(|\alpha | + |\beta |)/2}  | T_{I , h} \partial _x ^{\alpha }
 \partial _{\xi } ^{\beta }
 E_h ( X , S_I, T_I) | \leq K_3^{|I|} \sum _{(\gamma , \delta ) \in
 {\cal M}'_2 (I) }
 h^{(|\gamma | + |\delta |)/2}  |  \partial _x ^{\gamma } \partial
 _{\xi } ^{\delta }
 E_h ( X , S_I, T_I) |, \el
where ${\cal M}'_2 (I)=\{(\gamma , \delta)\in{\cal M}_2 (I)\ |\  \gamma _j + \delta _j \geq 2,\ \forall j\}$.
 $j$.
Thus, in view of  (\ref{5.17}), (\ref{5.18}), (\ref{5.13}) and (\ref{5.15}),
 $$ |T_{I , h} G_{I , h} (X)| \leq K_4^{|I|} \sum _{(\gamma , \delta
 ) \in {\cal M}'_2 (I) }
 h^{(|\gamma | + |\delta |)/2}  \Vert A_h\Vert
_{{\cal L}_{|\gamma  | + |\delta |}  ( Q_A) } \prod _{j\in I}
 Q_A(e_j , 0)^{\gamma _j /2}  Q_A(0 , e_j)^{\delta _j /2}. $$
If  $(\gamma , \delta )\in {\cal M}'_2 (I)$, one has $|\gamma |
 + |\delta |
 \leq 4 |I|\leq 4d$ and $2 \leq \gamma _j + \delta _j \leq 4$,  $\forall j\in I$.
Consequently, if $\lambda _j$ and $S$ are defined by (\ref{5.6}),
 $  Q_A(e_j , 0)^{\gamma _j /2}  Q_A(0 , e_j)^{\delta _j /2}  \leq
 \lambda _j S $.
Since $|{\cal M}_2 (I)|=9^{|I|}$,
one has with another universal constant   $K_5$, if
 $h <1$,
 $$|T_{I , h} G_{I , h} (X)| \leq  \Vert A_h\Vert
_{{\cal L}_{4d}  ( Q_A) } \prod _{j\in I} (K_5 h S \lambda _j ). $$
Summing over all subsets  $I\subset\{ 1,
\dots  , d \}$. One obtains,
$$ |F_{ h} (X)| \leq  \Vert A_h\Vert _{{\cal L}_{4d}  ( Q_A) }
 \prod _{j\leq d}  ( 1  +   h K S^2 \lambda _j).$$
Proceeding similarly for the derivatives, one deduces
(\ref{5.5}). Finally, in view of (\ref{2.16}) and (\ref{5.4}) the operators  $Op_h^{weyl}(F_h)$ and $ A_h $ having the same Wick symbol are actually equal, which proves (\ref{5.3}).

\hfill $\Box$

\begin{theo}\label{t1.8}  Let  $H$ be a real Hilbert space with finite dimension $d$, $A\in{\cal L}(H^2)$, $A^*=A$, $A\geq 0$, and $Q_A$ be the quadratic form (\ref{1.2}). Let $m\geq 0$
and $(A_h)$ a family of operators in ${\cal L} _{m+4d} (Q_A)
 $.
 Then, there exists $(F_h)$ belonging to $S_m (H^2 , Q_A)$  satisfying,
 $A_h = Op_h ^{weyl }(F_h)$ and,
\bl{1.11} \Vert F_h \Vert _{ m, Q_A} \leq \Vert A_h \Vert _{ {\cal L}
 _{m+4d} (Q_A) }
 \   D(H^2/{\rm Ker} A\ ,  \ q_A\  ,\   K  h \max (1, \Vert {\cal F}
 A \Vert _{H^2 , q_A } ) \ ,
 \ ({\cal F} A) )^{1/2},\el
with some universal constant $K>0$.
\end{theo}

{\it Proof of Theorem \ref{t1.8}.} \newline
 Let $H$ be a finite dimensional Hilbert space, $d={\rm dim}\,H$, and $A$ a linear map in $H^2$ being nonnegative and symmetric. Let $(e_j )$ be an orthonormal basis of $H$, let $u_j =
(e_j, 0)$ and $v_j = (0, e_j)$. Let $p\leq d$ and
$(U_j) , (V_j) $, $j\in\{1,\dots ,d\})$, a symplectic basis of $H^2$ satisfying properties in Theorem \ref{t3.1}. There exists a symplectic map  $\chi $  such that  $\chi (u_j) = U_j,\chi (v_j)= V_j$. Let $U_{\chi}$ be the metaplectic operator in
   $L^2 (H, \mu _{H, h/2})$
of Theorem 18.5.9 in \cite{Hor} satisfying
   (\ref{3.3}). From Proposition \ref{p4.3}, if the familly $(A_h)\in {\cal L} _{m+4d }
(Q)$ then the family  $D_h =  U_{\chi }^{\star} A_h U_{\chi }$ belongs to ${\cal L} _{m+4d } (Q \circ \chi )$. If $p$ is the integer of Theorem \ref{t3.1}, one has for $j\leq p$, 
$$ \Vert [ L_h (v_j), D_h] \Vert \leq  h \Vert  D_h \Vert _{m+4d, Q
\circ \chi }
Q \circ \chi (v_j)^{1/2}  = h \Vert  D_h \Vert _{m+4d, Q \circ \chi }
Q (V_j)^{1/2}= 0.$$
Consequently, $D_h$ commutes with $L_h(v_j)$ which is from (\ref{4.4}) a multiplication by the coordinate $x_j$. Set $x =(x', x'')$ with $x' = ( x_1 , \dots  ,x_p)$ and $x'' = (x_{p+1} , \dots  x_d)$. For any $x'\in\R^p$, there exists
 $D_h (x')$ in $L^2(\R^{d-p} , \mu _{\R^{d-p}  ,
h/2})$  satisfying,
$$ (D_hf )(x' , x'') = \Big ( D_h (x') f( x', \cdot ) \Big ) (x'').
$$
From Proposition \ref{p5.2} applied with $\R^{d-p}$ and with $D_h (x')$ together with the quadratic form $\widetilde{Q}$ being the restriction of $Q\circ\chi$ to $\R^{2(d-p)}$, there is a function on $\R^{2(d-p)}$
such that $D_h (x') = Op_h^{weyl } (K_{ h}(x'))$. Moreover,
$$ \Vert K_{ h}(x')\Vert _{m, Q\circ \chi|_{\R^{2(d-p)}}} \leq
\Vert D_h \Vert _{{\cal L}_{m+  4d }  ( Q_A \circ \chi ) }
 \prod _{p < j\leq d}  ( 1  +   h K S^2 \lambda _j),  $$ 
 where $K$ is a universal constant, and the $\lambda_j$ are defined as in $(\ref{5.6})$.
Set,
$$G_h (X' , X'')= G_h ( x' , \xi' , x'', \xi'') = K_h (x' ) ( x'' ,
\xi '').$$
This function is independant of  $\xi'$.  We have $Op_h ^{weyl}
(G_h) = D_h$. In addition,
$$ \Vert G_h \Vert _{m,  Q_A \circ \chi } \leq
\Vert D_h \Vert _{{\cal L}_{m+  4n }  ( Q_A \circ \chi ) }
 \prod _{p < j\leq n}  ( 1  +   h K S^2 \lambda _j).  $$
From (\ref{5.6}) applied to $Q_A \circ \chi$, one has
$$ \lambda _j = \max ( Q\circ \chi( e_j , 0), Q\circ \chi(0, e_j) )
= \max ( Q (U_j) , Q(V_j)).
$$
In view of Theorem \ref{t3.1}, if $j>p$, one has  $ Q_A(U_j)=
Q_A(V_j) = \lambda _j  $ where the $(\lambda _j )^2$ are the non vanishing eigenvalues of $ - ({\cal F} A)^2$. Consequently,
$$  \prod _{p < j\leq d}  ( 1  +   h K S^2 \lambda _j) =  D(H^2/{\rm
Ker} A\ ,  \ q_A\  ,\   K  h \max (1, \Vert {\cal F} A \Vert _{H^2 ,
q_A } ) \ ,
 \ ({\cal F} A) ).$$
In view of Proposition \ref{p4.3},
 $$ \Vert D_h \Vert _{{\cal L}_{m+  4d }  ( Q_A \circ \chi ) } =
 \Vert A_h \Vert _{{\cal L}_{m+  4d }  ( Q_A  ) }. $$
The function $F_h = G _h \circ \chi^{-1}  $ is in $S_m (H^2 , Q_A)$
and its norm verifies (\ref{1.11}). Since $Op_h ^{weyl} (G_h) = D_h$ and since  $D_h =  U_{\chi }^{\star} A_h U_{\chi }$, it comes  from
(\ref{3.3}) that $OP_h^{weyl}(F_h) = A_h$,
which proves Theorem \ref{t1.8}.

\hfill$\Box$

\subsection{The infinite dimensional case.}\label{s5.3}

 Let $(A_h)$ be a family of  ${\cal L} (Q_A) $.  Let $(e_j)$ be an Hilbert basis of $H$  such that  $e_j\in B'$ for every $j$ and let $E_n$ be the subspace generated by the $e_j$, $j\leq n$. The variable $X\in H^2$ may be written as
 $(X' , X'')$, with $X'\in E_n^2$, $X''\in (E_n^{\perp})^2$. Let $\widetilde E_n^{\perp}$ be the space of
 $x\in B$ such that $u(x) = 0$ for all $u \in E_n$.

{\it Step 1.}  By proposition \ref{p3.3},  $(i, E_n^{\perp} ,
\widetilde E_n^{\perp})$ is a Wiener space. For any $X''\in (E_n^{\perp})^2$ one may define a coherent state with $H$ replaced by $E_n^{\perp}$, $\Psi _{X'',
h}$ belonging to $L^2 (  \widetilde E_n^{\perp}  , \mu
_{ \widetilde E_n^{\perp}  , h/2})$. Let $A_{h ,E_n,  X''}$ the map defined by,
$$ <A_{h ,E_n,  X''} f, g> =  < A_h ( f \otimes \Psi _{X'', h} , g
\otimes \Psi _{X'', h}  >. $$
The operator $A_{h ,E_n,  X''}$ is bounded in
$L^2 ( E_n , \mu _{E_n , h/2})$. Let $Q_{A, E_n}=Q_A|_{E_n^2}$, let us prove that,
\bl{5.19} \Vert  A_{h , E_n, X''} \Vert _{{\cal L} ( Q_{A , E_n} ) } \leq
   \Vert  A_{h } \Vert _{{\cal L} ( Q_A) }.\el
If $V\in E_n^2$, if $f$ and $g$ are polynomial functions, one has,
 $$ < [L(V), A_{h , E_n,  X''} ] f, g > = < [L(V) , A_h ] ( f \otimes \Psi _{X'', h} ) \ , \ g \otimes \Psi _{X'', h} )> . $$
As $\Psi _{X'', h}$ is of norm 1, one has for any  $X''\in (E_n^{\perp})^2$,
  $$ \Vert [L(V), A_{h , E_n,  X''} ] \Vert \leq \Vert [L(V), A_{h }
  ] \Vert .$$
Considering in the same way the iterated commutators, one proves (\ref{5.19}). From Theorem \ref{t1.8} there exists $F_{ h , E_n , X''} (X')$ on $ E_n^2$  such that $ Op_h ^{weyl, E_n} ( F_{ h , E_n , X''} (X')) = A_{h , E_n,
  X''}$.
From Theorem \ref{t1.8}, one also has,
$$\Vert F_{ h, E_n, X'' } \Vert _{  Q_{A , E_n} } \leq \Vert  A_{h
 , E_n,X''}\Vert _{ {\cal L}
 (Q_{A, E_n}) }
  \ D_n^{1/2},  $$
with
 $D_n =  D(E_n^2/{\rm Ker} A_n\ ,  \ q_A\  ,\   K  h \max (1, \Vert
 {\cal F} A_n \Vert _{E^2 , q_A } ) \ ,
 \ ({\cal F} A_n) ), $
where $A_n = \pi_{E_n} Ai_{E_n}$. Setting $F_{n, h} (X' , X'')
 =  F_{ h , E_n , X''} (X')$, one has $ H_{E_n^2 , h/2}  F_{ h, E_n, X'' }(X') =  \sigma_h^{wick , E_n} (
 A_{h , E_n,X''} ) (X')$,
so that,
\bl{5.20} H_{E_n^2 , h/2} F_{n , h}  = \sigma_h^{wick } (A_h).\el
We also have,
\bl{est}  \Vert F_{n , h} \Vert _{\infty} \leq   \Vert  A_{h } \Vert
   _{{\cal L} ( Q) }
     D_n.    \el
For each $V_1$, \dots , $V_p\in H^2$, the symbol  $(h/i)^ {p}  (d^p F_{n , h}) (\cdot)  (V_1$, \dots , $V_p)$
is associated in the same way to the family of operators
${\rm ad} ( L_h V_1) \dots  {\rm ad} ( L_h V_p) A_h$. This relation between commutators
and differentials is common to the Wick and Weyl symbols. Therefore, applying (\ref{est}) to this
operator, we have,
$$  h^p  \Vert (d^p F_{n , h}) (\cdot)  (V_1, \dots , V_p)  \Vert _{\infty}
 \leq D_n \Vert {\rm ad} ( L_h V_1) \dots  {\rm ad} ( L_h V_p)  A_{h } \Vert
   _{{\cal L} ( Q) }. $$
 Thus, by Definitions \ref{d1.7} and \ref{d1.3},
\bl{5.21} \Vert  F_{n , h} \Vert _Q \leq  \Vert  A_{h } \Vert
   _{{\cal L} ( Q) }
     D_n. \el

  {\it Step 2.}
One may prove as in Theorem \ref{t1.6},
\bl{5.22} \lim _{n\rightarrow \infty } D_n   =
 D(B_A ,  q_A , K  h \max (1, \Vert {\cal F} A \Vert _{B_A , q_A } )
 \ ,   \ ({\cal F} A) ). \el
Let us prove that $(F_{n, h})$ is a Cauchy sequence in $C_b(H^2)$. If $m<n$, let $S_{m,n}$ be the
 orthogonal supplement of $E_m$ in $E_n$. Using (\ref{5.20}) and $ H_{E_n ^2, h/2}  = H_{S_{m,n}^2,
 h/2} H_{E_m^2 , h/2}   $ one deduces that $ H_{E_m^2 , h/2}   ( F_{m, h} - H_{S_{m,n}^2, h/2} F_{n, h}) =
 0 $. Since
  $H_{E_m^2 , h/2}$ is injective, one has
 $  F_{m, h} = H_{S_{mn}^2, h/2} F_{n, h}$.  Applying
 Proposition \ref{p3.5} one sees,
   $$  \lim _{m \rightarrow \infty , m<n} \Vert F_{m, h}  - F_{n, h}
   \Vert_{\infty}  =
    \lim _{m \rightarrow \infty , m<n} \Vert(H_{S_{m,n}^2, h/2} - I)
    F_{n, h} \Vert_{\infty}
  \leq  \lim _{m \rightarrow \infty , m<n}  (h/2)^{1/2} \Vert F_{n,
  h} \Vert_{Q}
  \ {\rm Tr} ( A \pi _ {S_{m,n}^2})^{1/2}.   $$
By (\ref{5.21}) and (\ref{5.22}),  $(F_{nh})$ is bounded in $S(H^2, Q)$, and the last trace tends to 0 by  Theorem 6.3 of Gohberg Krein \cite{GK}. Therefore, $F_{n, h}$ has a limit, denoted by $F_h$, in $C_b(H^2)$.
As in the first step, we may apply the same arguments with Beals commutators of $A_h$
and to differentials of $F_h$. Therefore
     $F_h\in  S (H^2, Q)$ and  estimate (\ref{1.12}) still holds true.

   {\it Step 3.} Let us prove that $OP_h^{weyl} (F_h) = A_h $.
   Using the partial heat operator defined in (\ref{2.15}) and (\ref{3.10}), we see that,
$H_{h/2} = H_{ (E_n^2)^{\perp}, h/2} H_{ E_n^2, h/2}$,    and (\ref{5.20}) proves
\bl{heat} H_{h/2} F_{n  h} =  H_{ (E_n^2)^{\perp}, h/2} \sigma_h^{wick}
   (A_h) .\el
For the left hand side, we have,
   $$ \lim _{n\rightarrow \infty} \Vert H_{h/2} (  F_{n  h} - F_h )
   \Vert _{\infty}
   \leq  \lim _{n\rightarrow \infty} \Vert   F_{n  h} - F_h  \Vert
   _{\infty}
   = 0. $$
For the right hand side of (\ref{heat}), we remark, by Proposition \ref{p5.1},  that $\sigma_h^{wick} (A_h)$
is in $S(H^2,Q)$, and that, by Proposition \ref{p3.5},
   $$  \Vert ( H_{ (E_n^2)^{\perp}, h/2} -
   I)  \sigma_h^{wick} (A_h)
   \Vert _{\infty}
   \leq (h/2)^{1/2} \Vert \sigma_h^{wick} (A_h) \Vert _{ Q_A }
   {\rm Tr} ( A \pi_{(E_n^2)^{\perp}} ) ^{1/2}  .  $$
The last trace goes to 0 in view of \cite{GK}
  (Theorem 6.3). Therefore, when $n$ goes to infinity in (\ref{heat}),
  we obtain  $ H_{h/2} F_{ h}=
   \sigma_h^{wick} (A_h)$. Since $F_{ h}\in S(H^2,Q)$, we can define the Weyl operator
   $Op_h^{weyl} (F_h)$. By (\ref{2.16}), and by the previous equality, its  Wick symbol is the same  as the Wick symbol of
   $A_h$. Therefore $Op_h^{weyl} (F_h) = A_h$. Theorem \ref{t1.9} is thus proved.

\hfill$\Box$

\section{Operators composition.}\label{s6}

If  $(A_h)$ and $(B_h)$ are two families in ${\cal L} (Q)$ one sees, for any vectors $V_1$, \dots , $V_m\in H^2$,
 $$ \Vert {\rm ad} ( L_h V_1) \dots  {\rm ad} ( L_h V_m) (A_h \circ B_h)
 \Vert  \leq
  \ 2^m h^m  \Vert A_h \Vert _{{\cal L}(Q)  }  \  \Vert B_h \Vert
  _{{\cal L}(Q)  }
  \ \prod _{j= 1}^m Q_A( V_j ) ^{1/2},\quad h\in (0, 1]. $$
This may be read as,
\bl{6.1} \Vert (A_h \circ B_h)   \Vert_{{\cal L}(4Q)  }  \leq  \Vert A_h
\Vert _{{\cal L}(Q)  }  \  \Vert B_h \Vert _{{\cal L}(Q)  }.\el

 {\it Proof of Theorem \ref{t1.11}.} \newline
Let  $(F_h)$ and $(G_h)$ two families in $S(H^2 , Q_A) $.
From  (\ref{1.9}), the families of operators
 $A_h = Op_h^{weyl} (F_h)$ and $B_h = Op_h^{weyl} (G_h)$  are in
 ${\cal L}  (Q)$.
 By (\ref{6.1}) and by Theorem   \ref{t1.9}, there exists a family $(K_h)$ in
 ${\cal S}  (H^2, 4Q) $  satisfying  (\ref{1.13}) and (\ref{1.15}).
The proof of (\ref{1.14}) is similar.

\hfill$\Box$

\section{Applications in field theory.}\label{s7}

\subsection{Hilbert spaces.}\label{s7.1}

In order to describe the free evolution of photons, without interaction, one uses a configuration Hilbert space consisting of maps $f \in L^2 (\R^3, \R^3)$  satisfying  $k\cdot f(k) =0$ almost everywhere. The Hilbert space ${\cal H}_{ph}$  describing the quantized field, with no interaction, at a given time is the symmetrized Fock space ${\cal F}_s ( H_{\bf C})$ associated with the above Hilbert space $H$. Thus, it also may be considered as an $L^2$ space, namely, $ L^2 ( B , \mu _{B ,h/2})$, where $(i, H, B)$ is a Wiener space. The space $H^2$ is named phase space.

The Hilbert space describing the states of  $N$ non interacting fixed particles with $1/2$ spin at a given time is  ${\cal H}_{sp}  =(
 \C^2  )^{\otimes N}$.

For the whole system (the quantized field and the $N$ particles in a constant magnetic field $\beta$), the Hilbert space is the tensor product ${\cal H}_{ph} \otimes{\cal H}_{sp} $.

\subsection{Operators.}\label{s7.2} The Hamiltonian $H(h)$ of the system is a selfadjoint operator in  ${\cal H}_{ph} \otimes {\cal H}_{sp}$ which may be written as,
\bl{7.1} H(h) = H_{ph} \otimes I + h H_{int},\el
where $H_{ph}$ is the photons Hamiltonian (operator
in ${\cal H}_{ph}$). In the framework of Fock space it is defined as $H_{ph} = h d\Gamma (M)$ where $M:H_{\bf C} \rightarrow H_{\bf C} $ is the multiplication by $\omega(k) = |k|$ and where $d\Gamma $ associates with operators in  $H_{\bf C}$, selfadjoint unbounded operators in the Fock space ${\cal F}_s ( H_{\bf C})$  (see Reed-Simon \cite{RSII}, Section X). By the isomorphism between the Fock space ${\cal F}_s ( H_{\bf C})$ and $L^2 ( B , \mu _{B ,h/2})$, this also  defines a selfadjoint unbounded operator in  $L^2 (B , \mu _{B ,h/2})$. With this operator, one may associate a Wick symbol which is the following quadratic form defined on a dense subspace of  $H^2$,
$$H_{ph} (q , p) = {1\over 2} \int _{\R^3} |k| \Big [ |q(k)|^2 +  |p(k)|^2 \Big ] dk. $$

The interaction operator $H_{int}$ is obtained through three operators $B_j (x)$ $j=1,\dots ,3$, $x\in\R^3$ which are the components of the magnetic field.
These operators are Weyl operators related to symbols which are linear continuous forms on $H^2$ written as,
\bl{7.2}   B_j (x,  q, p) =  (q, p)  \cdot   B_{jx}, \el
where $ B_{jx}$ belongs to $H^2$ and written as, when identifying $H^2$ with the complexified $H_{\C}$,
\bl{7.3} B_{jx}(k) = {i\chi(|k|)|k|^{1\over 2} \over (2\pi)^{3\over 2}}
e^{-i(   k\cdot x   )} {k\wedge e_j \over |k|},\quad
k\in\R^3\backslash\{0\},\el
where $\chi \in {\cal S} (\R)$.

The operators $E_j(x)$ are also associated, in a similar way, with a linear form on $H^2$ written as,
\bl{7.4} E_j (x, q, p) =  (q, p)  \cdot   J B_{jx}, \el
where $J: H^2 \rightarrow H^2$ is the helicity operator defined by,
\bl{7.5}J(q, p) (k) = \left ( {k\wedge q(k) \over |k| } , {k\wedge p(k)\over |k| } \right ),\quad k \in \R^3 \setminus \{ 0 \}. \el
Let $\sigma _j$, $j=1,2,3$, be the Pauli matrices,
\bl{7.6} \sigma_1 = \begin{pmatrix}
0 & 1 \\ 1 & 0   \end{pmatrix}, \quad
   \sigma_2 = \begin{pmatrix} 0 & -i \\ i & 0   \end{pmatrix},
\quad
  \sigma_3 = \begin{pmatrix}  1 & 0 \\ 0 & -1  \end{pmatrix}.\ee

For $\lambda \leq N$ and $m\leq 3$, we denote by
$\sigma_m^{[\lambda]}$ the operator in ${\cal H} _{sp}$ defined by,
\bl{7.7} \sigma_m^{[\lambda]} = I \otimes \cdots \sigma_m \cdots \otimes I,\el
where $\sigma_m$ is located at   $\lambda ^{th}$ position. Let
$\beta = (\beta_1 , \beta _2, \beta_3)$ be
the external constant magnetic field. Let
$x_{\lambda }$, $\lambda=1,\dots ,N$, be the position in $\R^3$ where are localized the fixed $1/2$ spin particles. Then, the interaction Hamiltonian is written as,
\bl{7.8} H_{int} = \sum _{\lambda =1}^N  \sum _{m=1}^3 (\beta_m + B_m(x_{\lambda})) \otimes
\sigma_m^{[\lambda]}. \el
This operator $H_{int}$ is initially defined on the space of stochastic extensions of polynomials functions with values in ${\cal H}_{sp}$.

Let us prove that it has a selfadjoint extension with the same domain as the free operator  $H_0 = H_{ph } \otimes I$.  It is standard. It is known that if
$\varphi\in D(H_{ph}^{1/2})$ and any $V\in H^2$  satisfying
$V(k) / \sqrt {|k|}\in H^2$, one has with the notations above Definition \ref{d1.7},
 $$\Vert L_h(V)  \varphi||\leq  C  (|V/\sqrt{|k|}|
   \Vert H_{ph}^{1/2}\varphi \Vert  + C  |V|   \Vert \varphi \Vert.
 $$
See \cite{BFS} Lemma I.6.  Consequently, for any $\varepsilon
 >0$, there exists $C_\varepsilon>0$ (depending on $V$ and $h$) satisfying,
 $$\Vert  L_h(V)\varphi \Vert \leq \varepsilon \Vert H_{ph}\varphi
 \Vert  +
  C_\varepsilon \Vert \varphi \Vert,\quad   \varphi\in D(H_{ph}).  $$
 Since the operators  $B_j(x_{\lambda})$ have this form, using that  $B_{j, x_{\lambda} }$ defined in (\ref{7.3}) belongs to $H^2$
after  division by $|k|^{1/2}$, one has for $\phi\in
D(H_{0})$ and $\varepsilon >0$,
 $$ \Vert H_{int} \phi \Vert \leq \varepsilon  \Vert H_{0}\phi \Vert
  + C_\varepsilon  \Vert \phi \Vert.  $$
Consequently, by Kato-Rellich Theorem, the operator
 $H(h)$ defined in (\ref{7.1}) and (\ref{7.8}) has an extension as a selfadjoint operator with domain $D(H(h)) =
D(H_0) = D(H_{ph}) \otimes {\cal H}_{sp}$. For any $x\in\R^3$,
the operators $B_j(x)$ and $E_j (x)$ $j=1,\dots ,3$ are bounded from $D(H(h))$ to ${\cal H} _{ph} \otimes {\cal H}_{sp}$.

\subsection{Quadratic form on  $H^2$.}\label{s7.3}
We define, for any $t\in\R$, a quadratic form $Q_t$ on the phase space $H^2$. This form involves the free evolution of the magnetic field. In the Fock framework (see Reed-Simon Vol.II) one writes
$H_{ph} = h d\Gamma (M)$, where $M$ is the multiplication by $\omega
(k) = |k|$, and consequently, $e^{i{t \over h} H_{ph} } = \Gamma (
\chi_t)$ where $\chi_t $ is the multiplication by $e^{it |k|}$ in
$H_{\bf C}$. If one prefers not identifying $H^2$ with $H_{\C}$ then  $\chi_t$ is both unitary and symplectic in $H^2$ and is defined by,
\bl{7.9} \chi_t (q , p) = (q_t , p_t),\quad
\left \{ \begin{matrix} q_t (k) = \cos (t \omega (k)) q(k) +  \sin (t\omega (k)) p(k)\\ \\
 p_t (k) = -  \sin (t \omega (k)) q(k) +  \cos (t \omega (k)) p(k)
 \end{matrix} \right ..\el
Any operator of the form $U_\chi=\Gamma(\chi)$ with $\chi:H_{\C}\rightarrow H_{\C}$ being both $C-$linear and unitary, is written as in (\ref{4.10}). In other words,   $e^{i{t \over h} H_{ph} } $ is a metaplectic operator (or a Bogoliubov transform) corresponding to  the unitary and symplectic map $\chi_t $.

If  $F$ is a linear form on $H^2$ one has as in (\ref{4.6}), for any $t\in \R$,
\bl{7.10}  e^{i{t \over h}   H_{ph} } Op_ h ^{weyl} (F)    e^{-i{t \over h}
H_{ph} }
= Op_ h ^{weyl} (F \circ \chi_{t}).\el
For any
$x\in \R^3$ and $t\in \R$, the operator describing the free evolution of the field, without interaction, is
  \bl{7.11} B_j^{free} (x, t) =  e^{i{t \over h} H_{ph} } B_j (x) e^{-i{t
\over h} H_{ph} }.\el
 From (\ref{7.9}), (\ref{7.10}) and (\ref{7.11}), this operator is through the Weyl calculus associated with a symbol which is a continuous linear form on $H^2$ written as,
 \bl{7.12} B_j^{free} (x, t, q, p) =  (q, p) \cdot    B_{jxt}, \el
where $ B_{jxt}$ is an element of  $H^2\simeq H_{\bf C}$ given by,
\bl{7.13} B_{jxt}(k) = {i\chi(|k|)|k|^{1\over 2} \over (2\pi)^{3\over 2}}
e^{i( t|k | - k\cdot x   )} {k\wedge e_j \over |k|},\qquad
k\in\R^3\backslash\{0\}.\el

For  $t\in \R$, one may define a nonnegative quadratic form $Q_t$ on $H^2$ by,
\bl{7.14} Q_t (q, p) = 3N  |t| \sum _{m=1 }^3 \sum _{\lambda =1}^N \int_0^t
|   B_{m}^{free} (x_{\lambda}, s, q, p)|^2   ds.\el

\subsection{Statement of the results.}\label{s7.4}

At first, we consider the operator,
\bl{7.15} U_h^{red} (t)  = \left [ e^{i{t \over h}   H_{ph} } \otimes I \right ]   e^{-i{t \over h}  H(h) }.\el
It is proved in \cite{ALN-QED-A} that $U_h^{red} (t)$ is associated through the Weyl calculus with a function  $U (t, \cdot , h)$ in the class studied in  \cite{AJN-JFA}, but to this end, one  assumes that $\chi$ vanishes in a neighborhood of 0. Our goal is here to drop this condition and to prove that  $U_h^{red} (t)$ is associated with a function  $U (t, \cdot , h)$ in $S (H^2 , Q_t)$.

 \begin{theo}\label{t7.1} For any $t\in \R$, the family $   U_h^{red} (t)$, $h\in (0,1]$ defined in (\ref{7.15})
 belongs to the class ${\cal L} (Q_t)$ (see Definition \ref{d1.7}) with the quadratic form defined in (\ref{7.14}), taking values in ${\cal L} ({\cal H}_{sp})$. Moreover,
\bl{7.16} \Vert  U_h^{red} (t) \Vert _{{\cal L} (Q_t, {\cal L} ({\cal H}_{sp}) )} = 1.\el
 There exists a function  $U(t, \cdot , h)$ in
 $S(H^2, Q_t, {\cal L} ( {\cal H} _{sp}))$  satisfying,
\bl{7.17}   U_h^{red} (t)= Op_h^{weyl} ( U(t, \cdot , h).\el
Moreover, with the notations of Theorem \ref{t1.9},
\bl{7.18} \Vert  U(t, \cdot , h) \Vert _{Q_t}  \leq
\   D(B_{A_t} ,  q_{A_t} , K  h \max (1, \Vert {\cal F} A_t \Vert
_{B_{A_t} , q_{A_t} } ) \ ,
   \ ({\cal F} A_t) )^{1/2},\el
 where $A_t$ satisfies  $Q_t (X) =
 (A_tX)\cdot X$ for all $X$ in $H^2$ and $K$ is a universal constant. Let us observe that $A_t$ is a trace class operator.
 \end{theo}

We shall describe the time evolution of some observables, still in the framework of interaction between  $N$ fixed spin $1/2$ particles and the quantized field. These observables will be the three components of the particles spin,  the three components of the electric and magnetic fields at each point $x\in\R^3$. The operators  $E_m (x)$ being the components of the electric field are defined in Section \ref{s7.2}. One sets, for  $\lambda \leq N$,
\bl{7.19}  S_j^{[\lambda ] } (t, h) =e^{  i {t\over h} H(h)}    (I \otimes
 \sigma _j^{[\lambda]})
 e^{ - i {t\over h} H(h)}\el
and for  $x\in \R^3$,
\bl{7.20} B_j (x , t, h) =e^{  i {t\over h} H(h)}    (B_j (x) \otimes I)
 e^{ - i {t\over h} H(h)},\el
\bl{7.21} E_j (x , t, h) =e^{  i {t\over h} H(h)}    (E_j (x) \otimes I)
 e^{ - i {t\over h} H(h)}.\el
 We will also use the free evolution of these observables (without interaction), i.e.,
   $ B_j ^{free} (x , t, h)$ defined in  (\ref{1.16}) or in (\ref{7.11}) and of $ E_j ^{free} (x , t, h)$ for the electric field.

The following theorem describes the full evolution of the fields and spin operators.
 \begin{theo}\label{t7.2}  For $t\in \R$,
 $x\in \R^3$,  $\lambda \leq N$,  $1\leq j\leq 3$, $h\in (0,1]$, the families of operators  $S_j^{[\lambda ] } (t ,h)$, $ B_j  (x , t,
 h)-  B_j ^{free} (x , t, h)$,
 $ E_j  (x , t,   h)-  E_j ^{free} (x , t, h) $,
belong to
 ${\cal L} (4Q_t)$ (see Definition \ref{d1.7}) with $Q_t$ being the quadratic form defined in (\ref{7.14}), with values in  ${\cal L} ({\cal H}_{sp})$ . Their norms are satisfying,
\bl{7.22} \Vert S_j^{[\lambda ] } (t ,h)\Vert _{{\cal L} (4Q_t, {\cal L}
 ({\cal H}_{sp}) )} = 1, \el
\bl{7.23}\Vert  B_j  (x , t,   h)-  B_j ^{free} (x , t, h)\Vert _{{\cal L}
(4Q_t, {\cal L} ({\cal H}_{sp}) )} \leq  h Q_t ({\cal F}  B_{jxt}  )
^{1/2}, \el
\bl{7.24} \Vert E_j  (x , t,   h)-  E_j ^{free} (x , t, h) \Vert _{{\cal L}
(4Q_t, {\cal L} ({\cal H}_{sp}) )}
 \leq h    Q_t ({\cal F} J B_{jxt}  ) ^{1/2}. \el

For each $x\in \R^3$ and any $t\in \R$, for every  $\lambda
 \leq N$ and $j\leq 3$, there exist functions   $  S_j^{[\lambda ] } (t, \cdot ,h)$, $B_j ^{res}
 (x , t, \cdot , h)$,
 $E_j ^{res} (x , t, \cdot , h)$ belonging to $S( H^2 , 4 Q_t , {\cal L} ( {\cal H}
 _{sp}))$ satisfying,
\bl{7.25} S_j^{[\lambda ] } (t ,h) = Op_h ^{weyl} (    S_j^{[\lambda ] }
(t, \cdot ,h) ),\el
\bl{7.26}B_j  (x , t,   h)=  B_j ^{free} (x , t, h) +  h Op_h ^{weyl} ( B_j
^{res} (x , t, \cdot , h) ),\el
\bl{7.27} E_j  (x , t,   h)=  E_j ^{free} (x , t, h) +  h Op_h ^{weyl} ( E_j
^{res} (x , t, \cdot , h) ).\el
\end{theo}

Let $N$ be the  photon number operator defined in Fock spaces framework by  $N = d\Gamma (I)$ (see Reed-Simon \cite{RSII}).  The operator $N$ can also be view as an unbounded selfadjoint  operator in $L^2(B , \mu _{B , h/2})$. The Wick symbol of this operator is,
$$ N(q , p) ={1\over 2h} ( |q|^2 + |p|^2). $$

The next result is concerned by the photon number evolution, namely,
\bl{7.28} N( t, h) =e^{  i {t\over h} H(h)}    (N \otimes I)
 e^{ - i {t\over h} H(h)}.\el

\begin{theo}\label{t7.3} One has,
\bl{7.29} {d\over dt} N(t , h) =  \sum _{\lambda =1}^N  \sum _{m=1}^3  X_{m
\lambda } (t)
\circ S_m^{[\lambda ] } (t, h),  \el
where
$S_m^{[\lambda ] } (t, h)$ is defined in (\ref{7.19}) and where, with the notations of Definition \ref{d1.7},
\bl{7.30}  X_{m \lambda } (t) = -  L_h ( \chi_t B_{m , x_{\lambda }}  )
\otimes I  +  Y_{m \lambda } (t),\el
where $Y_{m \lambda }  (t)   \in S( H^2 , 4 Q_t , {\cal L} (
 {\cal H} _{sp}))$ and satisfies,
\bl{7.31} \Vert Y _{m \lambda } (t)  \Vert _{{\cal L} (4 Q_t, {\cal L}
({\cal H}_{sp}) )}
\leq h Q_t ( \chi_t B_{m , x_{\lambda }}  ) ^{1/2}.\el
There exists a family of functions  $R_{m, \lambda , h , t }$
belonging to $ S(4 Q_t, {\cal L} ({\cal H}_{sp}) ) $   satisfying,  
\bl{7.32}  Y _{m \lambda } (t)  = h Op_h^{weyl} ( R_{m, \lambda , h , t }).\el
\end{theo}

\subsection{Proofs.}\label{s7.5}

{\it Proof of Theorem \ref{t7.1}.}\newline
For  $V\in H^2$, let
$L_h(V)$ be the operator introduced  above Definition \ref{d1.7} and set
$\widetilde L_h(V) = L_h(V) \otimes I$. One denotes by $H_{int}$
the operator defined in  (\ref{7.8}) and by $H_{int}^{free} (t)$
the corresponding free evolution operator defined in (\ref{1.16}).  One notes that,
\bl{7.33}H^{free}_{int}(t)
= \sum _{\lambda = 1}^N \sum_{j=1}^3 ( \beta _j + B_j^{free}
(a_{\lambda} , t)\otimes  \sigma _j ^{[\lambda]}.\el 
Differentiating (\ref{7.15}), one observes that for any  $f\in D( H(h))$,
$$ {d \over dt} U^{red}_h(t)  f = - i H^{free}_{int}(t) U^{red}_h(t)
f. $$
Similarly to \cite{ALN-QED-A}, one deduces that, for any $V\in H^2$, 
$$  {\rm ad} \widetilde L_h(V )U_h^{red} (t)  =i ^{-1} \int _0^t
U_h(t , s) [ \widetilde L_h(V) , H^{free}_{int}(s)]U_h(s, 0) ds,  $$
where we set $U_h(t , s) = U_h^{red} (t) (U_h^{red} (s)
)^{\star}$. Indeed, both hand sides are solutions to the same differential system, namely,
$$  X'(t) = - i H^{free}_{int}(t)X(t) - i  [ \widetilde L_h(V) ,
H^{free}_{int}(t)]
 U_h^{red} (t), $$
with the same initial data $X(0) = 0$. Iterating, one proves, for any   $V_1$, \dots  ,$V_m$ in $H^2$, 
\bl{7.34}  {\rm ad} \widetilde L_h(V_1) \dots   {\rm ad} \widetilde L_h(V_m
)U_h^{red} (t)  =
 (-i)^m  \sum _{\varphi \in S_m}  \int _{  \Delta _ m(t)}
    U(t, s_n)[ \widetilde L_h(V_{\varphi (m)}) ,
    H^{free}_{int}(s_m)]\dots \el
$$ \dots   U(s_n, s_{m-1} )[ \widetilde L_h(V_{\varphi (m-1)}) ,
H^{free}_{int}(s_{n-1})]
 \dots  U(s_2, s_1)[ \widetilde L_h(V_{\varphi (1)}) ,
 H^{free}_{int}(s_1)]
U(s_1, 0) ds_1 \dots ds_m,$$ 
where $S_m$ is the set of permutations $\varphi$ of $\{ 1, \dots  , m \}$ and where, for $t>0$,
$$ \Delta _ m(t) = \{ (s_1, \dots , s_m) \in
\R^n , 0 < s_1 < \dots  < s_m < t  \}.  $$ 
In view of (\ref{7.12}), (\ref{7.13})
and (\ref{7.33}), one has for any $V\in H^2$ and $t\in \R$,
\bl{7.35} \Vert [ \widetilde L_h(V) , H^{free}_{int}(t)]
\Vert \leq h N_t (V),\quad N_t (V) = \sum _{j=1}^3 \sum
_{\lambda = 1}^N |V \cdot B_{j, a_{\lambda }, t}|, \el
where $B_{j, x, t}$ is defined in (\ref{7.13}). Since $U(t, s)$
is unitary, from  (\ref{7.34}) and  (\ref{7.35}), one deduces
$$ \Vert  {\rm ad} \widetilde L_h(V_1) \dots   {\rm ad} \widetilde
L_h(V_m )U_h^{red} (t)
\Vert \leq h^m  \sum _{\varphi \in S_m}  \int _{  \Delta _ m(t)}
N_{s_1} (V _{\varphi (1)} )  \dots N_{s_m} (V _{\varphi (m)} ) ds_1
\dots ds_m $$ 
$$ = h^m \prod _{j=1}^m \int _0^t N_s(V_j) ds \leq h^m \prod
_{j=1}^m
 Q_t (V_j) ^{1/2}. $$
One then obtains the statement  (\ref{7.16}) of Theorem \ref{t7.1} and the second statement follows using  Theorem \ref{t1.9}.

\hfill$\Box$

{\it Proof of Theorem \ref{t7.2}. }\newline
The operator $I \otimes \sigma
_j^{[\lambda]}$ commutes with  $H_{ph} \otimes I$. Thus,
$$ S_j^{[\lambda ] } (t, h)=  U_h^{red } (t) ^{\star} (I \otimes
\sigma _j^{[\lambda]}) U_h^{red } (t).$$
The point (\ref{7.22}) is then obtained from Theorem \ref{t7.1} and inequality (\ref{6.1}). For the magnetic field, one has
$ B_j (x , t, h) =  U_h^{red } (t) ^{\star} (B_j^{free} (x , t)
\otimes I) U_h^{red } (t)$,
and consequently, since  $U_h^{red} (t)$ is a unitary map,
$$  B_j (x , t, h) - B_j^{free} (x , t) \otimes I  =   U_h ^{red}
(t))^{\star} \circ
[ B_j^{free} (x , t) ,  U_h ^{red} (t)].  $$ 
With Definition \ref{d1.7} notations, one has $ B_j^{free} (x , t)= L_h {\cal F},
B_{jxt}$ where $B_{jxt}$ belonging to $H^2$ is written in (\ref{7.13}) when identifying $H^2$ with $H_{\C}$. Since the family
 $ U_h ^{red} (t)$ belongs to ${\cal L} (Q_t, {\cal L} ({\cal H}_{sp}) )$, with norm 1, the commutator  $[ B_j^{free} (x , t) ,  U_h ^{red} (t)]$ is in this space and
$$ \Vert [ B_j^{free} (x , t) ,  U_h ^{red} (t)] \Vert _{{\cal L}
(Q_t, {\cal L} ({\cal H}_{sp}) )}
\leq   h Q_t ({\cal F}  B_{jxt}  ) ^{1/2}. $$ 
The statement (\ref{7.23}) is obtained in view of estimate (\ref{6.1}). For the electric field, the proof is the same. The other statements follow from Theorem \ref{t1.9}.

\hfill$\Box$

{\it Proof of Theorem \ref{t7.3}. }\newline
  From (\ref{7.28}) one sees,
$$ N'(t, h) = {i\over h} e^{{it \over h}  H(h)} [H(h), N\otimes I]
e^{-{it \over h} H(h) }.  $$
In the Fock spaces setting,  $N = d\Gamma (I)$ and
$H_{ph} = h d\Gamma (M)$. Since the operators $I$ and $M$
commutes, we have $[H_{ph}, N] = 0$. From (\ref{7.1})
and (\ref{7.8}), 
$$ N'(t, h) = i e^{{it \over h}  H(h)} [H_{int}, N\otimes I]  e^{-{it
\over h} H(h) }  =  \sum _{\lambda =1}^N  \sum _{m=1}^3  X_{m \lambda
} (t)
\circ S_m^{[\lambda ] } (t, h),$$ 
where $S_m^{[\lambda ] } (t,
h)$ is defined in (\ref{7.19}) and with
$$ X _{m \lambda } (t)  = i  e^{{it \over h}  H(h)}
( [B_m ( x_{\lambda}) , N]\otimes I )   e^{-{it \over h} H(h) }. $$
In view of  (\ref{7.15}), 
$$ X _{m \lambda } (t)  = i  U_h^{red} (t) ^{\star}
\left ( e^{i{t \over h}   H_{ph} }[B_m ( x_{\lambda}) , N]
 e^{-i{t \over h}   H_{ph} } \otimes I \right )
 U_h^{red} (t).   $$
Using  \cite{RSII} (Theorem X.41) and   \cite{DG} in the Fock space setting,
$$ [B_m ( x_{\lambda}) , N] = i L_h B_{m , x_{\lambda }}.$$
 From the metaplectic covariance result (\ref{4.7}), one gets,
$$  e^{i{t \over h}   H_{ph} } (L_h B_{m , x_{\lambda }})
   e^{-i{t \over h}   H_{ph} }  =  L_h ( \chi_t B_{m , x_{\lambda }}
   ).$$
Consequently,
$$ X _{m \lambda } (t)  = -   U_h^{red} (t)^{\star}
( L_h ( \chi_t B_{m , x_{\lambda }}  )  \otimes I )
    U_h^{red} (t).$$
Since $U_h^{red} (t)$ is a unitary map, one deduces 
$$  X _{m \lambda } (t)  = -  L_h ( \chi_t B_{m , x_{\lambda }}  )
\otimes I    + Y _{m \lambda } (t),\quad  Y _{m \lambda } (t)  =  -  U_h^{red} (t)^{\star} \circ [
L_h ( \chi_t B_{m , x_{\lambda }}  )  \otimes I  , \ U_h^{red} (t)].
$$
Knowing that $ U_h^{red} (t)$ belongs to
${\cal L} (Q_t, {\cal L} ({\cal H}_{sp}) )$, with norm 1, the commutator in the above right hand side is the same class and one sees that,
$$ \Vert [  L_h ( \chi_t B_{m , x_{\lambda }}  )  \otimes I  ,
\ U_h^{red} (t)]   \Vert _{{\cal L} (Q_t, {\cal L} ({\cal H}_{sp}) )}
\leq h Q_t ( \chi_t B_{m , x_{\lambda }}  ) ^{1/2}.$$ 
From
(\ref{6.1}),  the family $Y _{m \lambda } (t) $ belongs to ${\cal L} (4Q_t,
{\cal L} ({\cal H}_{sp}))$ and
$$ \Vert Y _{m \lambda } (t)  \Vert _{{\cal L} (4 Q_t, {\cal L}
({\cal H}_{sp}) )}
\leq h Q_t ( \chi_t B_{m , x_{\lambda }}  ) ^{1/2}.$$ 
The last statement in Theorem \ref{t7.3} comes  from Theorem \ref{t1.9}.

\hfill$\Box$

\medskip

laurent.amour@univ-reims.fr\newline
LMR EA 4535 and FR CNRS 3399, Universit\'e de Reims Champagne-Ardenne,
 Moulin de la Housse, BP 1039,
 51687 REIMS Cedex 2, France.

richard.lascar@imj-prg.fr\newline
Institut Mathématique de Jussieu UMR CNRS 7586,  Analyse Alg\'ebrique, 4 Place Jussieu, 75252 Paris. Boite 247.

jean.nourrigat@univ-reims.fr\newline
LMR EA 4535 and FR CNRS 3399, Universit\'e de Reims Champagne-Ardenne,
 Moulin de la Housse, BP 1039,
 51687 REIMS Cedex 2, France.


\begin{thebibliography}{99}


\bibitem{AJN-JFA} L. Amour, L. Jager, J. Nourrigat, {\it On bounded Weyl
pseudodifferential operators in Wiener spaces,} Journal of Functional
Analysis 269 (2015),
 pp. 2747-2812.


\bibitem{ALN-BEALS} L. Amour, R. Lascar, J. Nourrigat, {\it Beals
characterization of pseudodifferential
 operators in Wiener spaces}, Appl. Math. Res. Express (2016).




\bibitem{ALN-QED-A} L. Amour, R. Lascar, J. Nourrigat, {\it Weyl calculus in QED
I. The unitary group,}
 preprint, arXiv:1510.05293, october 2015.










\bibitem{BFS} V. Bach, J. Fr\"ohlich, I. M. Sigal, {\it Quantum
electrodynamics of confined nonrelativistic particles,} Adv. Math.
137 (1998), no. 2, 299-395.

\bibitem{Bea} R. Beals, {\it  Characterization of pseudodifferential operators
and applications,}  Duke Math. J. 44 (1977), no. 1, 45-57



\bibitem{Bon1} J.M. Bony,
{\it Caract\'erisation des opd.}  S\'eminaire EDP, X. Expos\'e n°23, 17pp, (1996-1997).




\bibitem{Bon2} J.M. Bony,
{\it Characterization of pseudo-differential operators,}
 Progress in non linear differential equations and their applications. Vol. 84. Birkha\"user, 21-34. (2013).



\bibitem{BC} J.M. Bony, J.Y. Chemin,
{\it Espaces fonctionnels associ\'es au calcul de Weyl-H\"ormander,}
 Bull. Soc. Math. France. {\bf 122}, n°1 77-118, (1994).


\bibitem{BMK} L. Boutet de Monvel, P. Kree, {\it Pseudodifferential
operators and Gevrey classes.}  Ann. Institut Fourier 17, (1967),
295-323.



\bibitem{BD} L. Bruneau,  J. Derezi\'nski, {\it  Bogoliubov Hamiltonians and
one-parameter groups of Bogoliubov transformations.}  J. Math. Phys.
48 (2007), no. 2, 022101, 24 pp.




\bibitem{CV} A.P. Calder\'on, R. Vaillancourt, {\it A class of bounded
pseudo-differential operators,} Proc. Nat. Acad. Sci. U.S.A, {\bf
69}, (1972), 1185-1187.



\bibitem{CR} M. Combescure, D. Robert, {\it Coherent states and
 applications in mathematical physics.}  Theoretical
 and Mathematical Physics. Springer, Dordrecht, 2012. xiv+415 pp.
 ISBN: 978-94-007-0195-3




\bibitem{DG}  J. Derezi\'nski, C.  G\'erard,  {\it Asymptotic completeness
in quantum field theory.
 Massive Pauli-Fierz Hamiltonians. }
Rev. Math. Phys. 11 (1999), no. 4, 383-450.






\bibitem{GK} I. C. Gohberg,  M. G. Krein, {\it Introduction to the theory of
linear nonselfadjoint operators}.  Translations of Mathematical
Monographs, Vol. 18 American Mathematical Society, Providence, R.I.
1969.


\bibitem{Gro} L. Gross, {\it Abstract Wiener spaces}, Proc. 5th Berkeley Sym.
Math. Stat. Prob, {\bf 2}, (1965), 31-42.





\bibitem{HSSS} F.Hiroshima, I. Sasaki, H. Spohn, A. Suzuki, {\it Enhanced
Binding in Quantum Field Theory}, Kyushu University COE Lecture Note
38, 2012. arXiv:1203.1136.

\bibitem{Hor} L. H\"ormander, {\it The analysis of linear partial differential
operators,} Volume III, Springer, 1985.



\bibitem{Jag} L. Jager, {\it Stochastic extensions of symbols in Wiener spaces
and heat operator}, arXiv:1607.02253, july 2016.



\bibitem{Jan}  S. Janson  {\it Gaussian Hilbert spaces,} Cambridge Tracts in
Maths, {\bf 129}, Cambridge Univ. Press (1997).

\bibitem{Kuo} H. H. Kuo, {\it Gaussian measures in Banach spaces.}  Lecture
Notes in Mathematics, Vol. 463. Springer, Berlin-New York, 1975.


\bibitem{Las} B. Lascar, {\it Une classe d'op\'erateurs
elliptiques du second ordre sur un espace de Hilbert,}  J. Funct.
Anal. {\bf 35} (1980), no. 3, 316-343.


\bibitem{Ler} N. Lerner, {\it Metrics on the phase space and non-selfadjoint pseudo-differential operators,} Pseudo-Differential Operators. Theory and Applications, 3. Birkhäuser Verlag, Basel, 2010.




\bibitem{Ram} R. Ramer, {\it On nonlinear
Transformations of Gaussian measures, } J. Funct. Analysis, 15
(1974), 166-187.



\bibitem{RSII} M. Reed, B.Simon, {\it Methods of modern mathematical physics,}
Academic Press, New York, London, 1978.




\bibitem{Sha} D. Shale. {\it Linear symmetries of free boson fields.}
Trans. Amer. Math. Soc., 103, 149-167, 1962.



 \bibitem{Sim-PP2} B. Simon,  {\it   The $P(\varphi)_2$ Euclidean (Quantum)
 Field Theory,}  Princeton Series in
Physics, Princeton University Press, 1974.

\bibitem{Tay}  M. E. Taylor, {\it Pseudodifferential operators.}  Princeton Mathematical Series,
34. Princeton University Press, Princeton, N.J., 1981.

\bibitem{Unt-OH} A. Unterberger, {\it Oscillateur harmonique et op\'erateurs
pseudo-diff\'erentiels,}
  Ann. Inst. Fourier (Grenoble) 29 (1979), no. 3, xi, 201-221.

\bibitem{Zwo} M. Zworski,  {\it Semiclassical analysis,}
Graduate Studies in Mathematics, 138. American Mathematical Society, Providence, RI, 2012.



\end{thebibliography}
\end{document}